\documentclass[12pt,a4paper]{article}
\usepackage{amssymb}
\usepackage{amsfonts}

\usepackage{amsmath}
\usepackage{amsthm}
\usepackage{enumerate}

\RequirePackage[OT1]{fontenc}
\RequirePackage{amsthm,amsmath, xcolor}
\RequirePackage[numbers]{natbib}
\RequirePackage[colorlinks,citecolor=blue,urlcolor=blue]{hyperref}
\usepackage{amscd,amsfonts,amssymb,latexsym,array,hhline,graphicx}
\usepackage{float}

\newtheorem{theorem}{Theorem}[section]
\newtheorem{proposition}[theorem]{Proposition}

\newtheorem{lemma}[theorem]{Lemma}

\newtheorem{remark}{Remark}[section]
\newtheorem{corollary}[theorem]{Corollary}

\newcommand\cA{{\cal A}}

\newcommand\cG{{\cal G}}
\newcommand\cF{{\cal F}}
\newcommand\cL{{\cal L}}
\newcommand\cB{{\cal B}}

\newcommand\cN{{\cal N}}

\newcommand\cX{{\cal X}}
\newcommand\cD{{\cal D}}
\newcommand\cQ{{\cal Q}}

\newcommand\cR{{\cal R}}

\newcommand\ve{\varepsilon}

\newcommand\Er{\mbox{Err}}

\def\bbr{{\mathbb R}}

\def\text#1{\hbox{#1}}
\def\proof{{\noindent \bf Proof. }}
\def\endproof{\mbox{\ $\qed$}}

\def\E{{\bf E}}

\def\P{{\bf P}}
\def\B{{\bf B}}
\def\p{{\bf p}}
\def\g{{\bf g}}

\def\C{{\bf C}}
\def\D{{\bf D}}

\def\G{{\bf G}}
\def\L{{\bf L}}
\def\U{{\bf U}}
\def\M{{\bf M}}

\def\c{{\bf c}}

\def\R{{\bf R}}

\def\r{{\bf r}}

\def\Chi{{\bf 1}}

\def\d{\mathrm{d}}
\def\build #1_#2{\mathrel{\mathop{\kern 0pt #1}\limits_{#2}}}
\newcommand\tr{\mbox{tr}}
\newcommand\Trg{\mbox{Tr}}
\newcommand{\wh}{\widehat}
\newcommand{\wt}{\widetilde}
\newcommand{\zs}[1]{{\mathchoice{#1}{#1}{\lower.25ex\hbox{$\scriptstyle#1$}}
{\lower0.25ex\hbox{$\scriptscriptstyle#1$}}}}

\numberwithin{equation}{section}

\begin{document}
\title{
Adaptive
model selection method
 for a conditionally
Gaussian semimartingale regression in
 continuous time
\thanks{
This work is supported
 by RSF, Grant no 17-11-01049.
}}
\author{Pchelintsev E.A.,
\thanks{
Department of Mathematical Analysis and Theory of Functions,
 Tomsk State University,
 e-mail: evgen-pch@yandex.ru
}
 \and
 Pergamenshchikov S.M.\thanks{
 Laboratoire de Math\'ematiques Raphael Salem,
   Universit\'e de Rouen
   and
International Laboratory of Statistics of Stochastic Processes and
Quantitative Finance of Tomsk State University,
 e-mail: Serge.Pergamenchtchikov@univ-rouen.fr}
}
 \date{}
\maketitle

\begin{abstract}
This paper considers the problem of robust adaptive efficient estimating of a
periodic function in a continuous time regression model with the
dependent noises given by a general square integrable
semimartingale with a conditionally Gaussian distribution. An example of such
noise is the non-Gaussian Ornstein--Uhlenbeck--L\'evy processes. An adaptive model selection
procedure, based on the improved weighted least square estimates, is proposed.
Under some conditions on the noise distribution, sharp oracle inequality for the robust
risk has been proved and the robust efficiency of the model selection procedure
has been established. The numerical analysis results are given.
\end{abstract}

\\
\par
{\bf UDC : 519.2}

\bibliographystyle{plain}

\newpage

\section{Introduction}\label{sec:In}

Consider a regression model in continuous time
 \begin{equation}\label{sec:In.1}
  \d y_t=S(t)\d t+\d \xi_t\,,
  \quad 0\le t\le n\,,
 \end{equation}
where $S$ is an unknown $1$-periodic $\bbr\to\bbr$ function,
$S\in\L_\zs{2}[0,1]$; $(\xi_\zs{t})_\zs{t\ge 0}$ is an
unobservable conditionally Gaussian semimartingale  with the values in the
Skorokhod space $\D[0,n]$
 such that, for any cadlag $[0,n]\to\bbr$
function $f$ from $\L_\zs{2}[0,n]$, the stochastic integral
\begin{equation}\label{sec:In.2}
I_\zs{n}(f)=\int^n_\zs{0}f(s)\d \xi_\zs{s}
 \end{equation} is
well defined and has the following properties
\begin{equation}\label{sec:In.3}
\E_\zs{Q} I_\zs{n}(f)=0 \quad\mbox{and}\quad \E_\zs{Q}
I^2_\zs{n}(f)\le \varkappa_\zs{Q} \int^{n}_\zs{0}\,f^{2}(s)\d s
\,.
\end{equation}
 Here $\E_\zs{Q}$ denotes the expectation with respect to the distribution
 $Q$ of the noise process $(\xi_\zs{t})_\zs{0\le t\le n}$ on the space $\cD[0,n]$;
 $\varkappa_\zs{Q}>0$ is
some positive constant depending on the distribution $Q$.
The noise distribution $Q$ is unknown and assumed to belong to some
 probability family $\cQ_\zs{n}$ specified below. All necessary
tools concerning the stochastic calculus can be found, for example, in \cite{JacodShiryaev2002}.

The class of the disturbances $\xi$ satisfying conditions
\eqref{sec:In.3} is rather wide and comprises, in particular, the
L\'evy processes which are used in different applied problems
(see \cite{Be}, for details). The models \eqref{sec:In.1}
with the L\'evy's type noise naturally arise
 in the nonparametric functional statistics problems
(see, for example, \cite{FeVi, KoPe2009a, KoPe2009b}). Moreover, as is shown in Section~\ref{sec:Ex},
non-Gaussian Ornstein--Uhlenbeck-based models,
introduced in \cite{BaNi}, enter this class.
It is well-known that in the filtration theory the assumption
of the conditional gaussinity of the unobserved process with respect to the
observed one led to the extension of the classical Kalman--Bucy problem with the
closed form solution to a class of the stochastic models described by the equations
including nonlinearly the process under observation \cite{LiptserShiryaev1977}.

The problem is to estimate the unknown function $S$ in the model
\eqref{sec:In.1} on the basis of observations
$(y_\zs{t})_\zs{0\le t\le n}$.

 We define the error of an estimate $\wh{S}$ (any real-valued function measurable with
respect to $\sigma\{y_\zs{t}\,,\,0\le t\le n\}$) for $S$ by its integral
quadratic risk
\begin{equation}\label{sec:In.4}
\cR_\zs{Q}(\wh{S},S):=
\E_\zs{Q,S}\,\|\wh{S}-S\|^2\,,
\end{equation}
where $\E_\zs{Q,S}$ stands for the expectation with respect to the
distribution $\P_\zs{Q,S}$ of the process \eqref{sec:In.1} with a fixed
distribution $Q$ of the noise $(\xi_\zs{t})_\zs{0\le t\le n}$
and a given function $S$; $\|\cdot\|$ is the norm in $\cL_\zs{2}[0,1]$, i.e.
\begin{equation*}\label{sec:In.5}
\|f\|^2:=\int^1_\zs{0}f^2(t) \d t\,.
\end{equation*}
Since in our case the noise distribution $Q$ is unknown, we will measure the quality of an estimate $\wh{S}$
by the robust risk defined as
\begin{equation}\label{sec:In.6}
\cR^{*}_\zs{n}(\wh{S},S)=\sup_\zs{Q\in\cQ_\zs{n}}\,
\cR_\zs{Q}(\wh{S},S)
\end{equation}
which assumes taking supremum of the error
\eqref{sec:In.4} over the whole family of admissible distributions $\cQ_\zs{n}$ (see, for example, \cite{GaPe2006}).

We will study the stated problem from the standpoint of the model selection approach.
approach. It will be noted that the
origin of this method goes back to papers by Akaike \cite{Ak} and Mallows \cite{Ma}.
The further progress has been made by Barron,
Birg\'e and Massart \cite{BaBiMa, Mas}, who developed a non-asymptotic model selection method which enables one to derive nonasymptotic oracle inequalities for nonparametric regression models with the i.i.d. Gaussian disturbances.
Fourdrinier and Pergamenshchikov \cite{FoPe} extended the
Barron–Birg\'e–Massart method to the models with the spherically symmetric dependent observations. The authors in \cite{KoPe2010} applied this method to the nonparametric problem of estimating a periodic function in a continuous time model with a Gaussian colored
noise. Unfortunately, the oracle inequalities obtained in these papers can not provide the efficient
estimation in the adaptive setting.
For constructing adaptive
procedures in our case one needs to use the approach based on the sharp oracle inequalities, proposed by
Galtchouk and Pergamenshchikov \cite{GaPe2009a, GaPe2009b} for the heteroscedastic
regression models in discrete time and which developed by Konev and Pergamenshchikov \cite{KoPe2012, KoPe2015}
for nonparametric regression models in continuous time.

The goal of this paper is to develop the adaptive robust efficient model selection method
for the regression \eqref{sec:In.1} with dependent noises having conditionally Gaussian distribution
using the improved estimation approach.
 This paper proposes the shrinkage least squares estimates which
enable us to improve the non-asymptotic estimation accuracy.
For the first time such idea was proposed by Fourdrinier and Pergamenshchikov
in \cite{FoPe} for regression models in discrete time and by Konev and
Pergamenshchikov in \cite{KoPe2010} for Gaussian regression models in continuous time.
We develop these methods for the general semimartingale regression models in continuous
time. It should be noted that for the conditionally Gaussian
regression models we can not use the well-known improved estimators proposed
in \cite{JamesStein1961} for Gaussian or spherically symmetric observations. To apply
the improved estimation methods to the non-Gaussian regression models in
continuous time one needs to use the modifications of the well-known James
- Stein estimators proposed in \cite{KPP2014, Pchelintsev2013} for parametric problems.
 We develop the new analytical tools
 which allow one to obtain the sharp non-asymptotic oracle inequalities
for robust risks  under general conditions on the distribution of the noise
 in the model \eqref{sec:In.1}. This method enables us to treat
 both the cases of dependent and independent observations from the same standpoint, it does not
assume the knowledge of the noise distribution and leads to the
efficient estimation procedure with respect to the risk
\eqref{sec:In.6}. The validity of the conditions, imposed on the
noise in the equation \eqref{sec:In.1} is verified for a non-Gaussian Ornstein--Uhlenbeck process.

The rest of the paper is organized as follows. In the next Section~\ref{sec:Ex}, we describe
the Ornstein--Uhlenbeck process as the example of a
semimartingale noise in the model \eqref{sec:In.1}. In Section~\ref{sec:Imp} we construct the shrinkage
weighted least squares estimates and study the improvement effect. In Section~\ref{sec:Mo} we construct the model
selection procedure on the basis
of improved weighted least squares estimates and state the main results in the form of oracle inequalities for the quadratic risk
\eqref{sec:In.4} and the robust risk \eqref{sec:In.6}.
In Section~\ref{sec:Ae} it is shown that the
proposed model selection procedure for estimating $S$ in \eqref{sec:In.1} is asymptotically efficient
with respect to the robust risk \eqref{sec:In.6}.
In Section~\ref{sec:Sim} we illustrate the performance of the proposed model selection procedure
through numerical simulations.
In Section~\ref{sec:Stc} we establish some properties of the stochastic integrals with respect to the non-Gaussian
Ornstein-–Uhlenbeck process \eqref{sec:Ex.1}. Section~\ref{sec:Prf} gives the proofs of the main results.
In the Appendix some auxiliary lemmas are given.


\section{Ornstein-Uhlenbeck-L\'evy process}\label{sec:Ex}

Now we consider the noise process $(\xi_\zs{t})_\zs{t\ge 0}$ in
\eqref{sec:In.1} defined by a non-Gaussian Ornstein--Uhlenbeck
process with the L\'evy subordinator. Such processes are used in
the financial Black--Scholes type markets with  jumps (see, for
example, \cite{DeKl} and the references therein). Let the noise
process in \eqref{sec:In.1} obeys the equation
\begin{equation}\label{sec:Ex.1}
\d\xi_\zs{t} = a\xi_\zs{t}\d t+\d u_\zs{t}\,,\quad \xi_\zs{0}=0\,,
\end{equation}
where
\begin{equation}\label{sec:Ex.0+0}
u_\zs{t} =
\varrho_\zs{1}\, w_\zs{t}+\varrho_\zs{2}\,z_\zs{t}
\quad\mbox{and}\quad
z_\zs{t}=\int_\zs{\bbr}\,x*(\mu-\wt{\mu})_\zs{t}\,.
\end{equation}
Here $(w_\zs{t})_\zs{t\ge 0}$ is
a standard Brownian motion, $\mu(\d s\,\d x)$ is the jump measure with the deterministic
compensator $\wt{\mu}(\d s\,\d x)=\d s\Pi(\d x)$, $\Pi(\cdot)$ is a L\'evy measure, i.e.  some positive measure on $\bbr_\zs{*}=\bbr\setminus \{0\}$, see, for example
\cite{ContTankov2004, JacodShiryaev2002}, such that
\begin{equation}\label{sec:Ex.1-00_mPi}
\Pi(x^{2})=1
\quad\mbox{and}\quad
\Pi(x^{8})
\,<\,\infty\,.
\end{equation}
We use the notation $\Pi(\vert x\vert^{m})=\int_\zs{\bbr_\zs{*}}\,\vert z\vert^{m}\,\Pi(\d z)$. Note that the L\'evy measure
 $\Pi(\bbr_\zs{*})$ could be equal to $+\infty$.

 We assume that the nuisance parameters  $a\le 0$, $\varrho_\zs{1}$
and $\varrho_\zs{2}$ satisfy the conditions
\begin{equation}\label{sec:Ex.01-1}
-a_\zs{max}\le  a\le 0\,,\quad
0< \underline{\varrho}\le \varrho^{2}_\zs{1}
\quad\mbox{and}\quad
\sigma_\zs{Q}=\varrho^{2}_\zs{1}+\varrho^{2}_\zs{2}\,
\le
\varsigma^{*}
\,,
\end{equation}
where
 the bounds
$a_\zs{max}$, $\underline{\varrho}$ and $\varsigma^{*}$ are functions of $n$, i.e.
$a_\zs{max}=a_\zs{max}(n)$,
$\underline{\varrho}=\varrho_\zs{n}$
and $\varsigma^{*}=\varsigma^{*}_\zs{n}$, such that for any $\check{\delta}>0$
\begin{equation}\label{sec:Ex.01-2}
\lim_\zs{n\to\infty}\,\frac{a_\zs{max}(n)}{n^{\epsilon}}=0
\,,\quad
\liminf_\zs{n\to\infty}\,n^{\epsilon}\,
\underline{\varrho}_\zs{n}
\,>0
\quad\mbox{and}\quad
\lim_\zs{n\to\infty}\,n^{-\epsilon}\,\varsigma^{*}_\zs{n}
=0
\,.
\end{equation}

We denote by  $\cQ_\zs{n}$  the family
of all distributions of process \eqref{sec:In.1} -- \eqref{sec:Ex.1} on the Skorokhod space
$\D[0,n]$ satisfying the conditions \eqref{sec:Ex.01-1} -- \eqref{sec:Ex.01-2}.

It should be noted that in view of Corollary
\ref{Co.sec:Stc.1} the condition \eqref{sec:In.3} for the process \eqref{sec:Ex.1} holds with
\begin{equation}\label{sec:Ex.4}
\varkappa_\zs{Q}=2\varrho_\zs{*}
\,.
\end{equation}
Note also that the process
\eqref{sec:Ex.1} is conditionally-Gaussian square integrated semimartingale with respect to
$\sigma$-algebra $\cG=\sigma\{z_\zs{t}\,,\,t\ge 0\}$ which is generated by jump process $(z_t)_{t\ge 0}$.

\bigskip

\section{Shrinkage estimates}\label{sec:Imp}

For estimating the unknown function $S$ in \eqref{sec:In.1} we will consider it's Fourier expansion.
Let $(\phi_\zs{j})_\zs{j\ge\, 1}$ be an orthonormal basis in $\L_\zs{2}[0,1]$.
We extend these functions  by the periodic way on $\bbr$, i.e.  $\phi_\zs{j}(t)$=$\phi_\zs{j}(t+1)$ for any $t\in\bbr$.

$\B_\zs{1}$) {\em Assume that the basis functions are uniformly bounded, i.e.
for some  constant $\phi_\zs{*}\ge 1$, which may be depend on $n$,
}
\begin{equation}\label{sec:In.3-00}
\sup_\zs{0\le j\le n}\,\sup_\zs{0\le t\le 1}\vert\phi_\zs{j}(t)\vert\,
\le\,
\phi_\zs{*}
<\infty\,.
\end{equation}

\bigskip

$\B_\zs{2}$) {\em Assume that there exist some $d_\zs{0}\ge 7$
and $\check{a}\ge 1$
  such that
\begin{equation}\label{sec:In.3-01}
\sup_\zs{d\ge d_\zs{0}}\,
\frac{1}{d}
\,
\int^{1}_\zs{0}\,
 \Phi^{*}_\zs{d}(v)\,
\d v
\le\,
\check{a}
\,,
\end{equation}
where $\Phi^{*}_\zs{d}(v)=\max_\zs{t\ge v}
\left\vert
\sum^{d}_\zs{j=1}\,\phi_\zs{j}(t)\,\phi_\zs{j}(t-v)
\right\vert$.
}

 For example, we can take
 the trigonometric basis    defined as $\Trg_\zs{1}\equiv 1$ and for $j\ge 2$
\begin{equation}\label{sec:In.5}
 \Trg_\zs{j}(x)= \sqrt 2
\left\{
\begin{array}{c}
\cos(\varpi_\zs{j}x)\, \quad\mbox{for even}\quad j \,;\\[4mm]
\sin(\varpi_\zs{j}x)\quad\mbox{for odd}\quad j\,,
\end{array}
\right.
\end{equation}
where the frequency $\varpi_\zs{j}=2\pi[j/2]$ and  $[x]$ denotes integer part of $x$.

In Lemma \ref{Le.sec:A.0} we shown that these functions satisfy the condition $\B_\zs{2}$) with
\begin{equation}\label{sec:In.5-00}
d_\zs{0}=\inf\{d\ge 7\,:\,5+\ln d\le \check{a} d\}
\quad\mbox{and}\quad
\check{a}=(1-e^{-a_\zs{max}})/(4a_\zs{max})\,.
\end{equation}

We write the Fourier expansion of the unknown function $S$ in the form
$$
S(t)=\sum_{j=1}^\infty \theta_\zs{j}\phi_j(t),
$$
where the corresponding Fourier coefficients
\begin{equation}\label{sec:Imp.2}
\theta_\zs{j}=(S,\phi_j)= \int^1_\zs{0}\,S(t)\,\phi_\zs{j}(t)\,\d t
\end{equation}
can be estimated as
\begin{equation}\label{sec:Imp.3}
\wh{\theta}_\zs{j,n}= \frac{1}{n}\int^n_\zs{0}\,\phi_j(t)\,\d
y_\zs{t}\,.
\end{equation}
We replace the differential $S(t)\d t$ by the stochastic observed differential $\d y_\zs{t}$.
In view of \eqref{sec:In.1}, one obtains
\begin{equation}\label{sec:Imp.4}
\wh{\theta}_\zs{j,n}=\theta_\zs{j}+\frac{1}{\sqrt{n}}\xi_\zs{j,n}\,,
\quad
 \xi_\zs{j,n}=\frac{1}{\sqrt{n}}
I_\zs{n}(\phi_\zs{j})
\end{equation}
where $I_\zs{n}(\phi_\zs{j})$ is given in \eqref{sec:In.2}.
As in \cite{KoPe2012} we define a class of weighted least squares estimates for $S(t)$ as
\begin{equation}\label{sec:Imp.5}
\wh{S}_\zs{\gamma}=\sum^{n}_\zs{j=1}\gamma(j)\wh{\theta}_\zs{j,n}\phi_\zs{j}\,,
\end{equation}
where the weights $\gamma=(\gamma(j))_{1\leq j\leq n}\in\bbr^{n}$ belong to some finite set $\Gamma$ from $[0,\,1]^n$ for which we set
\begin{equation}\label{sec:Imp.6}
\nu=\mbox{card}(\Gamma)
\quad\mbox{and}\quad
\vert \Gamma\vert_\zs{*}=\max_{\gamma\in \Gamma} \, \sum_{j=1}^n\gamma(j)\,,
\end{equation}
where $\mbox{card}(\Gamma)$ is the number of the vectors $\gamma$ in $\Gamma$.
In the sequel we assume that all vectors from $\Gamma$ satisfies the following condition.

\bigskip

$\D_\zs{1}$) {\sl Assume that for
for any vector $\gamma\in\Gamma$ there exists
 some fixed integer $d=d(\gamma)$
 such that their first $d$ components
equal to one, i.e. $\gamma(j)=1$ for $1\le j\le d$ for any $\gamma\in\Gamma$. }

\bigskip

$\D_\zs{2}$) {\sl There exists $n_\zs{0}\ge 1$ such that
for any $n\ge n_\zs{0}$ there exists a  $\sigma$ - field $\cG_\zs{n}$
for which
 the random vector
 $\wt{\xi}_\zs{d,n}=(\xi_\zs{j,n})_\zs{1\le j\le d}$
 is the $\cG_\zs{n}$-conditionally Gaussian in $\bbr^{d}$ with the covariance matrix
\begin{equation}\label{sec:Imp.6-1}
\G_\zs{n}=\left(
\E\,\xi_\zs{i,n}\,\xi_\zs{j,n}|\cG_\zs{n})
\right)_\zs{1\le i,j\le d}
\end{equation}
and for some nonrandom
constant $l^{*}_\zs{n}>0$
\begin{equation}\label{sec:Imp.6-1-0}
\inf_\zs{Q\in\cQ_\zs{n}}\quad
\left(
\tr \,\G_\zs{n}
-
\,
\lambda_\zs{max}(\G_\zs{n})
\right)
\geq l^{*}_\zs{n}\quad \mbox{a.s.}\,,
\end{equation}
where $\lambda_\zs{max}(A)$ is the maximal eigenvalue of the matrix $A$.
}

As it is shown in Proposition~\ref{Pr.sec:Stc.2} the condition $\D_\zs{2}$)
holds for the non-Gaussian Ornstein--Uhlenbeck-based model \eqref{sec:In.1} -- \eqref{sec:Ex.1}.

Further,
for the first $d$ Fourier coefficients in \eqref{sec:Imp.4}
 we will use
 the improved estimation method proposed
for parametric models
 in
\cite{Pchelintsev2013}. To this end we
set $\wt{\theta}_\zs{n}=(\wh{\theta}_\zs{j,n})_\zs{1\le j\le d}$.
In the sequel we will use the norm $\vert x\vert^{2}_\zs{d}=\sum^{d}_\zs{j=1}\,x^{2}_\zs{j}$
for any vector $x=(x_\zs{j})_\zs{1\le j\le d}$ from $\bbr^{d}$.
Now
we define the shrinkage estimators as
\begin{equation}\label{sec:Imp.12}
\theta^{*}_\zs{j,n}=
\left(
1
-
g(j)
\right)
\,\wh{\theta}_\zs{j,n}\,,
\end{equation}
where  $g(j)=(\c_\zs{n}/|\wt{\theta}_\zs{n}|_\zs{d}) \Chi_\zs{\{1\le j\le d\}}$,
$$
\c_\zs{n}=
\frac{l^{*}_\zs{n}}{\left(r^{*}_\zs{n}+\sqrt{d\varkappa_\zs{*}/n}\right)\,n}
\quad\mbox{and}\quad
\varkappa_\zs{*}=\sup_\zs{Q\in\cQ_\zs{n}}\,\varkappa_\zs{Q}
\,.
$$
The positive parameter
 $r^{*}_\zs{n}$  is such that
\begin{equation}\label{sec:Imp.12+1}
\lim_\zs{n\to\infty}\,r^{*}_\zs{n}\,=\infty
\quad\mbox{and}\quad
\lim_\zs{n\to\infty}\,
\frac{r^{*}_\zs{n}}{n^{\check{\delta}}}
\,=\,0
\end{equation}
for any $\check{\delta}>0$.

Now we introduce a class of shrinkage
weighted least squares estimates for $S$ as
\begin{equation}\label{sec:Imp.11}
S^{*}_\zs{\gamma}=\sum^{n}_\zs{j=1}\gamma(j)\theta^{*}_\zs{j,n}\phi_\zs{j}\,.
\end{equation}

We denote the difference of quadratic risks of the estimates \eqref{sec:Imp.5} and \eqref{sec:Imp.11} as
$$
\Delta_{Q}(S):=\cR_\zs{Q}(S^{*}_\zs{\gamma},S)-\cR_\zs{Q}(\wh{S}_\zs{\gamma},S)\,.
$$
For this difference we obtain the following result.

\begin{theorem}\label{Th.sec:Imp.1}
Assume that the conditions $\D_\zs{1})$ --  $\D_\zs{2})$ hold. Then for any $n\ge n_0$
\begin{equation}\label{sec:Imp.11+1}
\sup_{Q\in\cQ_\zs{n}}\,\sup_\zs{\Vert S\Vert\le r^{*}_\zs{n}}
\Delta_{Q}(S)<-\c^2_\zs{n}
\,.
\end{equation}
\end{theorem}

\begin{remark}\label{Re;sec:Imp.1}
The inequality \eqref{sec:Imp.11+1} means that non asymptotically, i.e. for any $n\ge n_\zs{0}$
 the estimate \eqref{sec:Imp.11}   outperforms in mean square accuracy the estimate \eqref{sec:Imp.5}.
\end{remark}

\section{Model selection method and oracle inequalities}\label{sec:Mo}

This Section gives the construction of a model selection procedure  for
estimating a function $S$ in \eqref{sec:In.1} on the basis of improved weighted least square estimates and states
the sharp oracle inequality for the robust risk of proposed procedure.

\noindent The model selection procedure for the unknown function
$S$ in \eqref{sec:In.1} will be constructed on the basis of
a family of  estimates $(S^{*}_\zs{\gamma})_\zs{\gamma\in\Gamma}$.

The performance of any estimate $S^{*}_\zs{\gamma}$ will be measured by the
empirical squared error
$$
\Er_\zs{n}(\gamma)=\|S^*_\zs{\gamma}-S\|^2.
$$
In order to obtain a good estimate, we have to write a rule to choose a weight vector
$\gamma\in \Gamma$ in \eqref{sec:Imp.11}. It is obvious, that the best way is to minimise
the empirical squared error with respect to $\gamma$. Making use the estimate definition
\eqref{sec:Imp.11} and the Fourier transformation of $S$ implies
\begin{equation}\label{sec:Mo.1}
\Er_\zs{n}(\gamma)\,=\,
\sum^{n}_\zs{j=1}\,\gamma^2(j)(\theta^*_\zs{j,n})^2\,-
2\,\sum^{n}_\zs{j=1}\,\gamma(j)\theta^*_\zs{j,n}\,\theta_\zs{j}\,+\,
\sum^{n}_\zs{j=1}\theta^2_\zs{j}\,.
\end{equation}
Since the Fourier coefficients $(\theta_\zs{j})_\zs{j\ge 1}$ are
unknown, the weight coefficients $(\gamma_\zs{j})_\zs{j\ge 1}$ can
not be found by minimizing this quantity. To circumvent this
difficulty one needs to replace  the terms
$\theta^*_\zs{j,n}\,\theta_\zs{j}$ by their estimators
$\wt{\theta}_\zs{j,n}$. We set
\begin{equation}\label{sec:Mo.2}
\wt{\theta}_\zs{j,n}=
\theta^*_\zs{j,n}\,\wh{\theta}_\zs{j,n}-\frac{\wh{\sigma}_\zs{n}}{n}
\end{equation}
where $\wh{\sigma}_\zs{n}$ is the estimate for the limiting variance
of $\E_\zs{Q}\,\xi^{2}_\zs{j,n}$ which we choose in the following form
\begin{equation}\label{sec:Mo.3}
\wh{\sigma}_\zs{n}=\sum_{j=[\sqrt{n}]+1}^n \wh{t}_\zs{j,n}^2\,,\quad
\wh{t}_\zs{j,n}=\int_0^1 \Trg_\zs{j}(t)\d y_t.
\end{equation}
For this change in the empirical squared error, one has to pay
some penalty. Thus, one comes to the cost function of the form
\begin{equation}\label{sec:Mo.4}
J_\zs{n}(\gamma)\,=\,\sum^{n}_\zs{j=1}\,\gamma^2(j)(\theta^*_\zs{j,n})^2\,-
2\,\sum^{n}_\zs{j=1}\,\gamma(j)\,\wt{\theta}_\zs{j,n}\,
+\,\rho\,\wh{P}_\zs{n}(\gamma)
\end{equation}
where $\rho$ is some positive constant,
$\wh{P}_\zs{n}(\gamma)$ is the penalty term defined as
\begin{equation}\label{sec:Mo.5}
\wh{P}_\zs{n}(\gamma)=\frac{\wh{\sigma}_\zs{n}\,|\gamma|^2_\zs{n}}{n}
\,.
\end{equation}

\noindent
Substituting the weight coefficients, minimizing the cost function
\begin{equation}\label{sec:Mo.6}
\gamma^*=\mbox{argmin}_\zs{\gamma\in\Gamma}\,J_n(\gamma)\,,
\end{equation}
in \eqref{sec:Imp.5} leads to the improved model selection procedure
\begin{equation}\label{sec:Mo.7}
S^*=S^*_\zs{\gamma^*}\,.
\end{equation}
It will be noted that $\gamma^*$ exists because
 $\Gamma$ is a finite set. If the
minimizing sequence in \eqref{sec:Mo.6} $\gamma^*$ is not
unique, one can take any minimizer.

To prove the sharp oracle inequality, the following conditions will be needed
for the family $\cQ_n$ of distributions of the
noise $(\xi_\zs{t})_\zs{t\ge 0}$ in \eqref{sec:In.1}.

We need to impose some stability conditions for the noise Fourier transform sequence
$(\xi_\zs{j,n})_\zs{1\le j\le n}$ introduced in \cite{PchPerSISP2018}. To this end   for some parameter $\sigma_\zs{Q}> 0$
we set
the following function
\begin{equation}
\label{L_1_Q}
\L_\zs{1,n}(Q)=
 \,
\sum^{n}_\zs{j=1}\,
\left|
\E_\zs{Q}\,\xi^{2}_\zs{j,n}
-
\sigma_\zs{Q}
\right|
\,.
\end{equation}
In \cite{KoPe2012}
the  parameter $\sigma_\zs{Q}$ is called proxy variance.

\noindent $\C_\zs{1})$ {\it There exists a proxy variance
$\sigma_\zs{Q}> 0$
such that for any $\epsilon>0$
$$
\lim_\zs{n\to\infty}\frac{\L_\zs{1,n}(Q)}{n^{\epsilon}}
=0\,.
$$
}

\noindent Moreover, we define
\begin{equation*}
\label{L_2_Q}
\L_\zs{2,n}(Q)=
 \sup_\zs{|x|\le 1}
\E_\zs{Q}\,
\left(
\sum^{n}_\zs{j=1}\,x_\zs{j}\,
\wt{\xi}_\zs{j,n}
\right)^2
\quad\mbox{and}\quad
\wt{\xi}_\zs{j,n}
=\xi^2_\zs{j,n}-\E_\zs{Q} \xi^2_\zs{j,n}
\,.
\end{equation*}
\noindent $\C_\zs{2})$ {\it Assume that for any $\epsilon>0$
$$
\lim_\zs{n\to\infty}\frac{\L_\zs{2,n}(Q)}{n^{\epsilon}}
=0\,.
$$
}

\vspace{2mm}

As is shown in Propositions \ref{Pr.sec:L1-10} and \ref{Pr.sec:L2-11},
 both conditions $\C_\zs{1})$
and $\C_\zs{2})$ hold for the model \eqref{sec:In.1} with Ornstein-Uhlenbeck noise process \eqref{sec:Ex.1}.


\begin{theorem}\label{sec:Mo.Th.1}
If the conditions $\C_\zs{1})$ and $\C_\zs{2})$ hold for the
distribution $Q$ of the process $\xi$ in \eqref{sec:In.1}, then, for any $n\geq1$ and $0<\rho<1/2$,
the risk \eqref{sec:In.4} of estimate \eqref{sec:Mo.7} for $S$
satisfies the oracle inequality
\begin{equation}\label{sec:Mo.15_OrIneq}
\cR_\zs{Q}(S^{*},S)\,\le\, \frac{1+5\rho}{1-\rho}
\min_\zs{\gamma\in\Gamma} \cR_\zs{Q}(S^*_\zs{\gamma},S)
+
\frac{\B_\zs{n}(Q)}{\rho n}
\,,
\end{equation}
where
$\B_\zs{n}(Q)=\U_\zs{n}(Q)
\left(
1+\vert\Gamma\vert_\zs{*}\E_\zs{Q}|\wh{\sigma}_\zs{n}-\sigma_\zs{Q}|
\right)
$
and
 the coefficient $\U_\zs{n}(Q)$ is such that for any $\epsilon>0$
\begin{equation}
\label{termB_rest}
\lim_\zs{n\to\infty}\frac{\U_\zs{n}(Q)}{n^{\epsilon}}=0
\,.
\end{equation}
\end{theorem}

\bigskip

\noindent
 In the case, when the value of $\sigma_\zs{Q}$ in $\C_\zs{1})$ is known, one can take
$\wh{\sigma}_\zs{n}=\sigma_\zs{Q}$ and
\begin{equation}\label{sec:Mo.9}
P_\zs{n}(\gamma)=\frac{\sigma_\zs{Q}\,|\gamma|^2_\zs{n}}{n}
\end{equation}
and then we can rewrite the oracle inequality \eqref{sec:Mo.15_OrIneq}  with $\B_\zs{n}(Q)=\U_\zs{n}(Q)$.

\noindent
Now we study the estimate \eqref{sec:Mo.3}.

\begin{proposition}\label{sec:Mo.Prop.1}
Let in the model \eqref{sec:In.1} the function $S(\cdot)$ is continuously differentiable.
Then, for any $n\geq 2$,
\begin{equation*}\label{sec:Mo.15_kappa}
\E_\zs{Q}|\wh{\sigma}_\zs{n}-\sigma_\zs{Q}| \leq
\frac{\varkappa_\zs{n}(Q)(1+\|\dot{S}\|^2)}{\sqrt{n}}\,,
\end{equation*}
where the term $\varkappa_\zs{n}(Q)$ possesses  the property  \eqref{termB_rest}
 and $\dot{S}$ is the derivative of the function $S$.
\end{proposition}
\noindent
To obtain the oracle inequality for the robust risk \eqref{sec:In.6}
 we need some additional condition on the distribution family $\cQ_\zs{n}$.
We set
\begin{equation}
\label{sigma*_n}
\varsigma^{*}=\varsigma^{*}_\zs{n}=\sup_\zs{Q\in\cQ_\zs{n}}
\sigma_\zs{Q}
\quad\mbox{and}\quad
\L^{*}_\zs{n}=
\sup_\zs{Q\in\cQ_\zs{n}}(\L_\zs{1,n}(Q)+\L_\zs{2,n}(Q))
\,.
\end{equation}

\noindent $\C^{*}_\zs{1}$) {\it Assume that the conditions $\C_\zs{1})$--$\C_\zs{2})$
hold and for any  $\epsilon>0$
$$
\lim_\zs{n\to\infty}\frac{\L^{*}_\zs{n}+\varsigma^{*}_\zs{n}}{n^{\epsilon}}
=0\,.
$$
}

\noindent
Now we impose the conditions on the set of the weight coefficients $\Gamma$.

\noindent $\C^{*}_\zs{2})$ {\it Assume that the set $\Gamma$ is such that for any $\epsilon>0$
\begin{equation*}\label{sec:Mo.8+1}
\lim_\zs{n\to\infty}\frac{\nu}{n^{\epsilon}}=0
\quad\mbox{and}\quad
\lim_\zs{n\to\infty}\,\frac{\vert\Gamma\vert_\zs{*}}{n^{1/2+\epsilon}}
=0\,.
\end{equation*}
}

\begin{theorem}\label{Th.sec:2.3}
Assume that the conditions
$\C^*_\zs{1})$--$\C^*_\zs{2})$ hold.  Then the robust risk
 \eqref{sec:In.6} of the estimate \eqref{sec:Mo.7} for
continuously differentiable function $S(t)$ satisfies for any $n\ge 2$  and
$0<\rho<1/2$ the oracle inequality
\begin{align*}\label{sec:Mo.20}
\cR^{*}_\zs{n}(S^{*},S)\,\le\,
\frac{1+5\rho}{1-\rho} \min_\zs{\gamma\in\Gamma}
\cR_n^{*}(S^*_\zs{\gamma},S) +\frac{1}{\rho n}\,\B^{*}_\zs{n}(1+\|\dot{S}\|^2)\,,
\end{align*}
where the term $\B^{*}_\zs{n}$ satisfies the property
\eqref{termB_rest}.
\end{theorem}

Now we specify the weight coefficients $(\gamma(j))_\zs{j\ge
1}$ in the way proposed in \cite{GaPe2009a}
 for a heteroscedastic regression
model in discrete time. Firstly, we define the normalizing coefficient which defined the minimax convergence rate
\begin{equation}
\label{upsilon-nnn}
 v_\zs{n}=
\frac{n}{\varsigma^{*}}
\,,
\end{equation}
where the upper proxy variance $\varsigma^{*}$ is defined in \eqref{sigma*_n}.  Consider a numerical grid of the form
\begin{equation*}\label{sec:Imp.7}
\cA_\zs{n}=\{1,\ldots,k^*\}\times\{r_1,\ldots,r_m\}\,,
\end{equation*}
where  $r_i=i\varepsilon$ and $m=[1/\varepsilon^2]$. Both
parameters $k^*\ge 1$ and $0<\varepsilon\le 1$ are assumed to be
functions of $n$, i.e. $ k^*=k^*(n)$ and
$\varepsilon=\varepsilon(n)$, such that for any $\delta>0$
\begin{equation*}\label{sec:Imp.8}
\left\{
\begin{array}{ll}
&\lim_\zs{n\to\infty}\,k^*(n)=+\infty\,,
\quad
\lim_\zs{n\to\infty}\,\dfrac{k^*(n)}{\ln n}=0\,,\\[6mm]
&\lim_\zs{n\to\infty}\varepsilon(n)=0
\quad\mbox{and}\quad
\lim_\zs{n\to\infty}\,n^{\delta}\varepsilon(n)\,=+\infty\, .
\end{array}
\right.
\end{equation*}
One can take, for
example,
$$
\varepsilon(n)=\frac{1}{\ln (n+1)}
\quad\mbox{and}\quad
k^*(n)=\sqrt{\ln (n+1)}\,.
$$
For each $\alpha=(\beta,r)\in\cA_\zs{n}$ we introduce the weight
sequence $\gamma_\zs{\alpha}=(\gamma_\zs{\alpha}(j))_\zs{j\ge 1}$
as
\begin{equation}\label{sec:Imp.9}
\gamma_\zs{\alpha}(j)=\Chi_\zs{\{1\le j\le d\}}+
\left(1-(j/\omega_\alpha)^\beta\right)\, \Chi_\zs{\{ d<j\le
\omega_\alpha\}}
\end{equation}
where $d=d(\alpha)=\left[\omega_\zs{\alpha}/\ln (n+1)\right]$,
$
\omega_\zs{\alpha}=\left(\tau_\zs{\beta}\,r\,v_\zs{n}\right)^{1/(2\beta+1)}$ and
$$
\tau_\zs{\beta}=\frac{(\beta+1)(2\beta+1)}{\pi^{2\beta}\beta}\,.
$$
We set
\begin{equation}\label{sec:Imp.10_Lambda}
\Gamma\,=\,\{\gamma_\zs{\alpha}\,,\,\alpha\in\cA_\zs{n}\}\,.
\end{equation}
It will be noted that such weight coefficients satisfy the condition $\D_\zs{1})$ and in this case the cardinal of the set $\Gamma$ is
$\nu=k^{*} m$. Moreover,
taking into account that $\tau_\zs{\beta}<1$ for $\beta\ge 1$
we obtain for the set \eqref{sec:Imp.10_Lambda}
\begin{equation*}
\label{sec:Ga.1++1--2}
 \vert \Gamma\vert_\zs{*}\,
 \le\,1+
\sup_\zs{\alpha\in\cA}
  \omega_\zs{\alpha}
\le 1+(\upsilon_\zs{n}/\ve )^{1/3}\,.
\end{equation*}

\begin{remark}
Note that the form \eqref{sec:Imp.9}
for the weight coefficients
 was proposed by Pinsker in \cite{Pi}
 for the efficient estimation in the nonadaptive case, i.e. when the regularity parameters of the function $S$ are known.
In the adaptive case  these weight coefficients are
  used in \cite{KoPe2012, KoPe2015}
   to show the asymptotic efficiency for model selection procedures.
\end{remark}

\bigskip
\section{Asymptotic efficiency}\label{sec:Ae}

In order to study the asymptotic efficiency we define the following functional Sobolev ball
\begin{equation}\label{sec:Ae.1}
W_\zs{k,\r}=\{f\in\C^{k}_\zs{p}[0,1]\,:\,
\sum_\zs{j=0}^k\,\|f^{(j)}\|^2\le \r\}\,,
 \end{equation}
where $\r>0$ and $k\ge 1$ are
some unknown parameters, $\C^{k}_\zs{p}[0,1]$ is the space of
 $k$ times differentiable $1$ - periodic $\bbr\to\bbr$ functions
 such that for any $0\le i \le k-1$
$$
f^{(i)}(0)=f^{(i)}(1)
\,.
$$
 In order to formulate our asymptotic results
we define the well-known Pinsker constant which gives the lower bound for normalized asymptotic risks
\begin{equation}\label{sec:Ae.3}
l_\zs{k}(\r)\,=\,((1+2k)\r)^{1/(2k+1)}\,
\left(\frac{k}{\pi (k+1)}\right)^{2k/(2k+1)}
\,.
\end{equation}

It is well known that for any $S\in W_\zs{k,\r}$
 the optimal rate of convergence is
$T^{-2k/(2k+1)}$ (see, for example, \cite{GaPe2009b}).
On the basis of the model selection procedure
we construct the adaptive
procedure $\wh{S}_\zs{*}$ for which we obtain the following asymptotic upper bound
for the quadratic risk.

\noindent Now we show that the parameter \eqref{sec:Ae.3} gives a lower bound
for the asymptotic normalized risks.
To this end we denote by $\Sigma_\zs{n}$ of all estimators $\wh{S}_\zs{n}$ of $S$ measurable with respect to
the process
\eqref{sec:In.1}, i.e.
$\sigma\{y_\zs{t}\,,\,0\le t\le n\}$.

\begin{theorem}\label{Th.sec: Ae.2}
The robust risk \eqref{sec:In.6}
admits the following lower bound
 \begin{equation}\label{sec:Ae.5}
\liminf_\zs{n\to\infty}\,
\inf_\zs{\wh{S}_\zs{n}\in\Sigma_\zs{n}}
\,v_\zs{n}^{2k/(2k+1)}
\sup_\zs{S\in W_\zs{k,\r}}\,\cR^{*}_\zs{n}(\wh{S}_\zs{n},S)
\,
\ge l_\zs{k}(\r) \,.
 \end{equation}
\end{theorem}

\noindent
We show that this lower bound is sharp in the following sense.

\begin{theorem}\label{Th.sec: Ae.1}
The quadratic risk \eqref{sec:In.2} for the
 estimating procedure $S^{*}$ has the following asymptotic upper bound
 \begin{equation}\label{sec:Ae.4}
\limsup_\zs{n\to\infty}\,v_\zs{n}^{2k/(2k+1)}
\sup_\zs{S\in W_\zs{k,\r}}\,\cR^{*}_\zs{n}(S^{*},S)
\,
\le l_\zs{k}(\r)
\,.
 \end{equation}
\end{theorem}

\noindent
It is clear that Theorem \ref{Th.sec: Ae.1}
 and Theorem \ref{Th.sec: Ae.2}
imply
\begin{corollary}\label{Co.sec: Ae.1}
The model selection procedure $S^{*}$
is efficient, i.e.
\begin{equation}\label{sec:Ae.5}
\lim_{n\to\infty}\,(v_\zs{n})^{\frac{2k}{2k+1}}\,
\sup_\zs{S\in W_\zs{k,\r}}\,\cR^{*}_\zs{n}(S^{*},S)\,
= l_\zs{k}(\r)
\,.
\end{equation}
\end{corollary}

\begin{remark}\label{Re.2.2}
Note that the equality
\eqref{sec:Ae.5} implies that the parameter
\eqref{sec:Ae.3} is the Pinsker constant
in this case (cf. \cite{Pi}).
\end{remark}

\begin{remark}\label{Re.2.3}
It should be noted that the equality \eqref{sec:Ae.5}
means that the robust efficiency holds with
the  convergence rate
$$
(v_\zs{n})^{\frac{2k}{2k+1}}
\,.
$$
It is well known that for the simple risks
 the optimal (minimax) estimation convergence rate
for the functions from the set  $W_\zs{k,\r}$
 is $n^{2k/(2k+1)}$ (see, for example, \cite{Pi},  \cite{Nu}).
So,
 if  the upper bound of the distribution variance $\varsigma^{*}\to 0$  as $n\to\infty$
 we obtain the more rapid rate, and
 if $\varsigma^{*}\to \infty$  as $n\to\infty$
we obtain the more slow rate. In the case when $\varsigma^{*}$ is constant the robust rate is the same as the classical non-robust convergence rate.
\end{remark}

\section{Monte Carlo simulations}\label{sec:Sim}

In this section we give the results of numerical simulations to assess the performance and improvement of
the proposed model selection procedure \eqref{sec:Mo.6}.
We simulate the model \eqref{sec:In.1} with
$1$-periodic function $S$ of the form
\begin{equation}\label{sec:Sim_Sign_11}
S(t)=t\,\sin(2\pi t)+t^2(1-t)\cos(4\pi t)
\end{equation}
on $[0,\,1]$ and the L\'evy noise process $\xi_\zs{t}$ is defined as
$$
\d \xi_\zs{t}=-\xi_\zs{t} \d t + 0.5\,\d w_\zs{t}+0.5\, \d z_\zs{t}\,.
$$
Here $z_\zs{t}$ is a compound Poisson process with intensity $\lambda=\Pi(x^2)=1$ and a Gaussian
$\cN(0,\,1)$ sequence $(Y_\zs{j})_\zs{j\ge1}$ (see, for example, \cite{KoPe2015}).

We use the model selection procedure   \eqref{sec:Mo.6} with the weights \eqref{sec:Imp.9} in which
$k^*=100+\sqrt{\ln (n+1)}$, $r_i=i/\ln (n+1)$, $m=[\ln^2 (n+1)]$, $\varsigma^*=0.5$ and $\delta=(3+\ln n)^{-2}$.
We define the empirical risk as
$$
\cR(S^*,\,S)=\frac{1}{p}\sum_\zs{j=1}^p \wh{\E}\left(S_n^*(t_j)-S(t_j)\right)^2\,,
$$
$$
 \wh{\E}\left(S_n^*(\cdot)-S(\cdot)\right)^2=
\frac{1}{N}\sum_\zs{l=1}^N \left(S_{n,l}^*(\cdot)-S(\cdot)\right)^2\,,
$$
where the observation frequency $p=100001$ and the expectations was taken as an
average over $N = 1000$ replications.


\begin{table}[h]
\label{Tab1}
\caption{The sample quadratic risks for different optimal $\gamma$}
\begin{center}
\begin{tabular}{|c|c|c|c|c|}
  \hline
   $n$              & 100 & 200 & 500 & 1000 \\ \hline
  $\cR(S^*_\zs{\gamma^*},\,S)$    & 0.0289 & 0.0089 & 0.0021 & 0.0011 \\ \hline
  $\cR(\wh{S}_\zs{\wh{\gamma}},\,S)$ & 0.0457 & 0.0216 & 0.0133 & 0.0098 \\ \hline
  $\cR(\wh{S}_\zs{\wh{\gamma}},\,S)/\cR(S^*_\zs{\gamma^*},\,S)$ & 1.6 & 2.4 & 6.3 & 8.9 \\
  \hline
\end{tabular}
\end{center}
\end{table}

\begin{table}[h]
\label{Tab2}
\caption{The sample quadratic risks for the same optimal $\wh{\gamma}$}
\begin{center}
\begin{tabular}{|c|c|c|c|c|}
  \hline
   $n$              & 100 & 200 & 500 & 1000 \\ \hline
  $\cR(S^*_\zs{\wh{\gamma}},\,S)$    & 0.0391 & 0.0159 & 0.0098 & 0.0066 \\ \hline
  $\cR(\wh{S}_\zs{\wh{\gamma}},\,S)$ & 0.0457 & 0.0216 & 0.0133 & 0.0098 \\ \hline
  $\cR(\wh{S}_\zs{\wh{\gamma}},\,S)/\cR(S^*_\zs{\wh{\gamma}},\,S)$ & 1.2 & 1.4 & 1.3 & 1.5 \\
  \hline
\end{tabular}
\end{center}
\end{table}

Table 1 gives the values for the sample risks of the improved estimate \eqref{sec:Mo.6}
and the model selection procedure based on the weighted LSE (3.15) from \cite{KoPe2012} for different numbers
of observation period $n$. Table 2 gives the values for the sample risks of the the model selection procedure based on the weighted LSE (3.15) from \cite{KoPe2012} and it's improved version for different numbers
of observation period $n$.

\begin{figure}[h!]
\label{fig3}
\centering
    \includegraphics[width=0.5\textwidth]{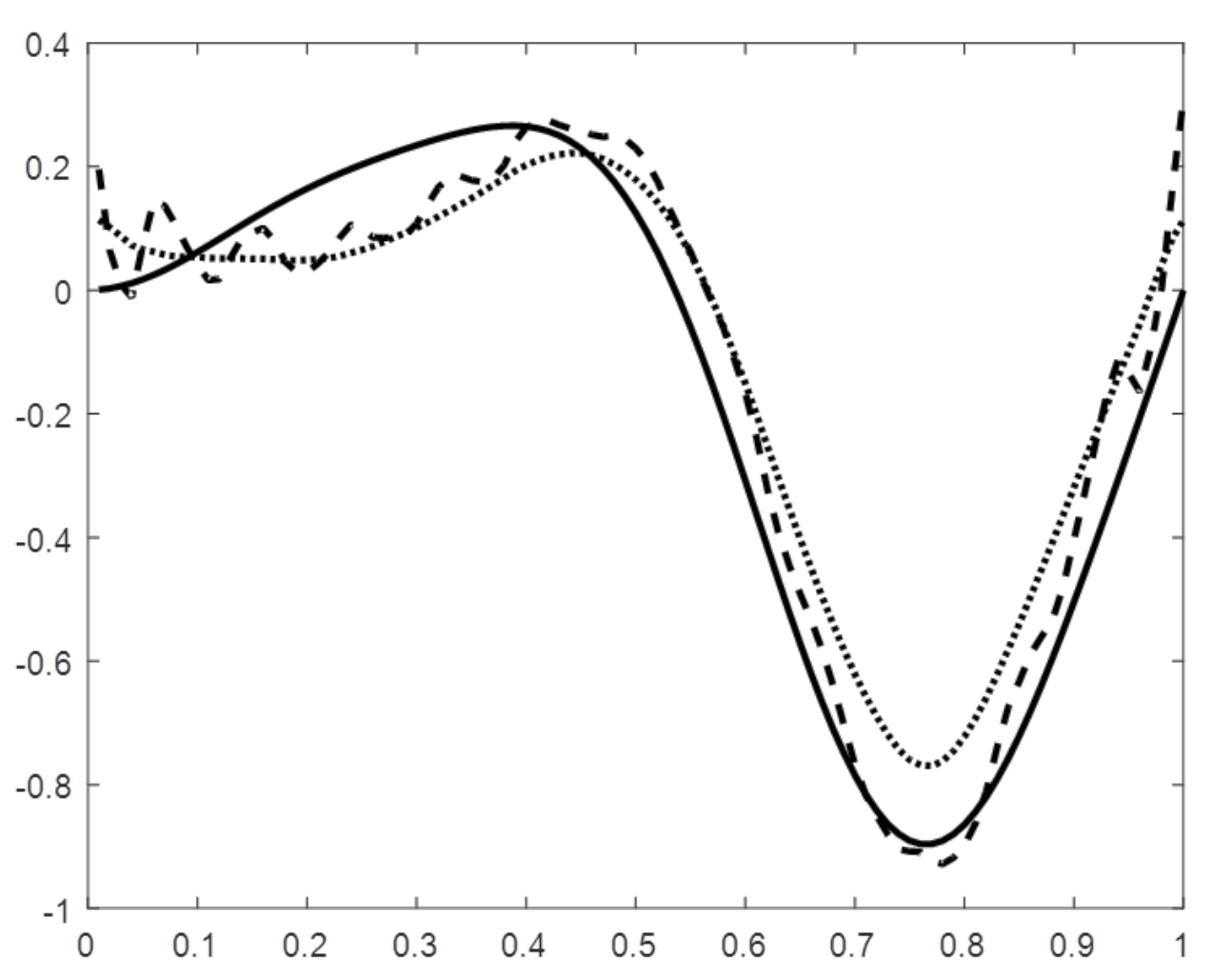}
  \caption{Behavior of the regression function and its estimates for $n=500$.}
\end{figure}


\begin{figure}[h!]\label{fig4}
\centering
 \includegraphics[width=0.5\textwidth]{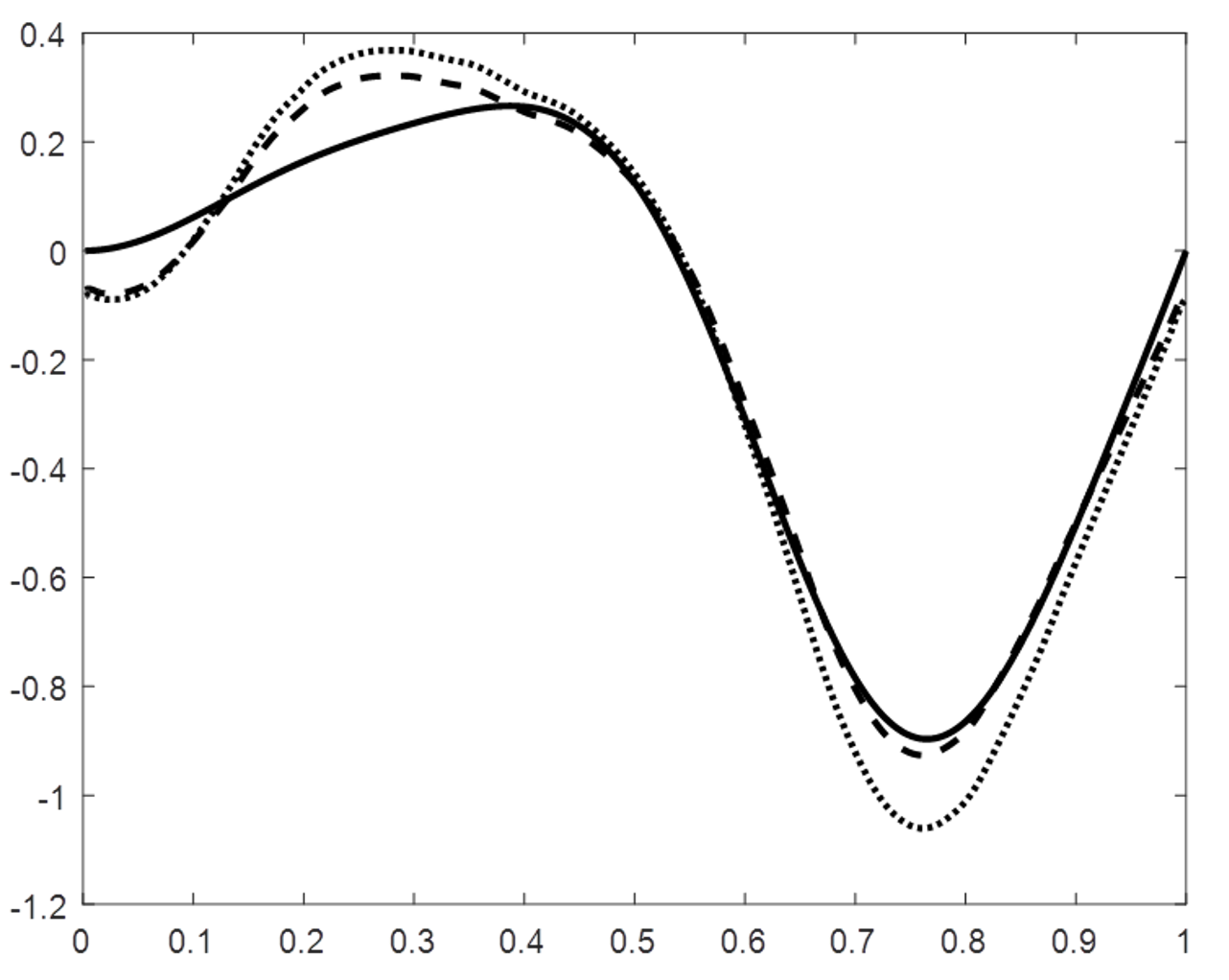}
  \caption{Behavior of the regression function and its estimates for $n=1000$.}
\end{figure}

\bigskip

\begin{remark}
Figures 1--2 show the behavior of
the procedures  \eqref{sec:Imp.5} and \eqref{sec:Mo.6}   depending
on the values of observation periods $n$.
The bold line is the function \eqref{sec:Sim_Sign_11},
the continuous  line is the  model selection
procedure based on the
least squares estimators $\wh{S}$
 and the dashed line is the improved model selection procedure $S^*$.
 From the Table 2 for the same $\gamma$ with various observations
numbers $n$ we can conclude that
theoretical result on the improvement effect \eqref{sec:Imp.11+1} is confirmed
by the numerical simulations.
Moreover, for the proposed shrinkage procedure, Table 1 and Figures 1--2,
we can conclude that the benefit is considerable for non large  $n$.
\end{remark}


\bigskip
\section{Stochastic calculus for Ornstein-Uhlenbeck-L\'evy process}\label{sec:Stc}

In this section we study the process \eqref{sec:Ex.1}.
\begin{proposition}\label{Pr.sec:Stc.1}
Let $f$ and $g$ be two nonrandom left continuous $\bbr_\zs{+}\to\bbr$ functions
with the finite right limits.
Then for any $t>0$
\begin{equation}\label{sec:Stc.1}
\E\, I_\zs{t}(f)I_\zs{t}(g)=\sigma_\zs{Q}\,
\tau_\zs{t}(f,g)\,,
\end{equation}
where
$\tau_\zs{t}(f,g)=
\int^{t}_\zs{0}\left(
f(s)g(s)
+
\check{\varepsilon}_\zs{s}(f)
g(s)+
f(s)
\check{\varepsilon}_\zs{s}(g)
\right)
\,\d s$
and
$$
\check{\varepsilon}_\zs{t}(f)=a\int^{t}_\zs{0}\,e^{a(t-s)}\,f(s)\left(
\frac{1+e^{2as}}{2}
\right)\,\d s
\,.
$$
\end{proposition}
\proof
Taking into account the definitions
\eqref{sec:Imp.4} and  \eqref{sec:Ex.1}
we obtain through the Ito formula that
\begin{align}\label{sec:Stc.2}
I_\zs{t}(f)\,I_\zs{t}(g)=
\sigma_\zs{Q}\,
\int^{t}_\zs{0}\,f(s)g(s)\d s+
a\int^{t}_\zs{0}
\Upsilon_\zs{s}(f,g)\,\xi_\zs{t}
\d s
+\M_\zs{t}(f,g)
\,,
\end{align}
where $\Upsilon_\zs{s}(f,g)=f(s) I_\zs{s}(g) + g(s) I_\zs{s}(f)$,
$$
\M_\zs{t}(f,g)=\int^{t}_\zs{0}\,\Upsilon_\zs{s-}(f,g)
\,
\d u_\zs{s}
+
\varrho^{2}_\zs{2}\,
\int^{t}_\zs{0}
f(s)\,g(s)\,
\d m_\zs{s}
$$
and $m_\zs{t}=x^{2}*(\mu-\wt{\mu})_\zs{t}$.
Moreover, using the Ito formula we obtain
\begin{equation}\label{sec:Stc.3}
\E\,I^{2}_\zs{t}(1)=
\E\, \xi^{2}_\zs{t}
=\sigma_\zs{Q}\,
\frac{e^{2at}-1}{2a}\,.
\end{equation}
Note now, that
\begin{align*}
\E\,
I^{2}_\zs{t}&=
\E
\left(a\int^{t}_\zs{0}\,f(s)\xi_\zs{s}\d s
+
\int^{t}_\zs{0}\,f(s)\,\d u_\zs{s}
\right)^{2}\\
&\le
2a^{2}\,
\int^{t}_\zs{0}\,f^{2}(s)\d s
\int^{t}_\zs{0}\,\E\,\xi^{2}_\zs{s}\d s
+2\sigma_\zs{Q}\,
\int^{t}_\zs{0}\,f^{2}(s)\d s\,.
\end{align*}
So, from here
\begin{equation}\label{sec:Stc.4}
\sup_\zs{0\le t\le n}
\E\,I^{2}_\zs{t}(f)
\le
2\sigma_\zs{Q}\,(\vert a\vert+
1)\,
\int^{n}_\zs{0}\,f^{2}(s)\d s
<\infty
\,.
\end{equation}
This implies immediately that
$\E\,\M_\zs{t}(f,g)=0$. Using this in \eqref{sec:Stc.2} yields
\begin{align}\nonumber
\E\,I_\zs{t}(f)\,I_\zs{t}(g)&=
\sigma_\zs{Q}\,
\int^{t}_\zs{0}\,f(s)g(s)\d s
\\[2mm]\label{sec:Stc.5}
&+
a\int^{t}_\zs{0}
\left(
f(s)\,\E \zeta_\zs{s}(g)
+
g(s)\,
\E \zeta_\zs{s}(f)
\right)
\,
\d s
\,,
\end{align}
where $\zeta_\zs{t}(f)=\xi_\zs{t}\,I_\zs{t}(f)=\,I_\zs{t}(1)\,I_\zs{t}(f)$. Therefore, putting $g=1$ in \eqref{sec:Stc.5}, we obtain that
$$
\E\,\zeta_\zs{t}(f)
=
\sigma_\zs{Q}\,
\int^{t}_\zs{0}\,f(s)\d s+
a\int^{t}_\zs{0}
\left(
f(s)
\E\,\zeta_\zs{s}(1)
+
\E\,\zeta_\zs{s}(f)
\right)
\,
\d s\,.
$$
Taking into account here, that $\zeta_\zs{t}(1)=\xi^{2}_\zs{t}$, we obtain that
$$
\E\,\zeta_\zs{t}(f)=
\sigma_\zs{Q}\,
\int^{t}\,e^{a(t-s)}f(s)
\frac{1+e^{2as}}{2}
\d s
=\sigma_\zs{Q}\,
\check{\varepsilon}_\zs{t}(f)\,.
$$
Therefore, using this in
 \eqref{sec:Stc.5} we obtain \eqref{sec:Stc.1}.
\endproof

\begin{corollary}\label{Co.sec:Stc.1}
For any cadlag function $f$  from $\L_\zs{2}[0,n]$
\begin{equation}\label{sec:Stc.6}
\E\, I^{2}_\zs{n}(f)
\le
2\sigma_\zs{Q}\,
\int^{n}_\zs{0}\,f^{2}(s)
\,\d s\,.
\end{equation}
\end{corollary}
\proof
Indeed, putting $f=g$ in \eqref{sec:Stc.1} we get
$$
\E\, I^{2}_\zs{n}(f)=\sigma_\zs{Q}\,
\int^{n}_\zs{0}\,\left(
f^{2}(t)
+2
\check{\varepsilon}_\zs{t}(f)\,
f(t)\right)
\,\d t\,.
$$
Moreover, note that
$$
\int^{n}_\zs{0}\,
\check{\varepsilon}_\zs{t}(f)\,
f(t)
\,\d t
=a\int^{n}_\zs{0}\,e^{ax}\,\int^{n}_\zs{x}\,\left(f(t)f(t-x)\frac{1+e^{2a(t-x)}}{2}\right)\d t
\,\d x\,.
$$
By the Bunyakovskii--Cauchy--Schwarz inequality
$$
\int^{n}_\zs{0}\,
\check{\varepsilon}_\zs{t}(f)
\,
\vert f(t)\vert
\,\d t
\le
\vert a\vert
\int^{n}_\zs{0}\,e^{ax}\,\,\d x\,
\int^{n}_\zs{0}\,f^{2}(t)\,\d t
\le
\int^{n}_\zs{0}\,f^{2}(t)\,\d t\,.
$$
This implies immediately upper bound
\eqref{sec:Stc.6}. Hence
Corollary \ref{Co.sec:Stc.1}.
\endproof

\noindent Now we set
\begin{equation}\label{sec:Stc.6-00}
\wt{I}_\zs{t}(f)=I^{2}_\zs{t}(f)
-\E\,I^{2}_\zs{t}(f)
\quad\mbox{and}\quad
V_\zs{t}(f)=\zeta_\zs{t}(f)\,
-
\E\,\zeta_\zs{t}(f)
\,.
\end{equation}
Using \eqref{sec:Stc.2} with $f=g$ we can obtain that
\begin{equation}\label{sec:Stc.6-01}
\d \wt{I}_\zs{t}(f)
=2af(t)\,V_\zs{t}(f)\d t
+
\d\wt{M}_\zs{t}(f)
\,,
\end{equation}
where
$\wt{M}_\zs{t}(f)=M_\zs{t}(f,f)$.
To study this process we need to introduce the following functions

\begin{equation}\label{sec:Stc.6-02-0}
\check{\tau}_\zs{t}(f,g)=f(t)g(t)\tau_\zs{t}(1,1)
+
f(t)\tau_\zs{t}(1,g)
+
g(t)\tau_\zs{t}(1,f)
+
\tau_\zs{t}(f,g)
\end{equation}
and
\begin{equation}\label{sec:Stc.6-01-00-1}
A_\zs{t}(f)=
\int^{t}_\zs{0} e^{3a(t-s)}f(s)\upsilon(s)\d s
+2
\,
\sigma_\zs{Q}^2
\,
\int^{t}_\zs{0} e^{3a(t-s)}\check{\varepsilon}_\zs{s}(f)\d s
\,,
\end{equation}
where $\upsilon(s)=a^{2}\E\,\wt{\xi}^{2}_\zs{s}+\sigma_\zs{Q}^2
\left( e^{2as}-1 \right) +a\check{\varrho}_\zs{2}$,
 $\wt{\xi}_\zs{s}=\xi^{2}_\zs{s}-\E\,\xi^{2}_\zs{s}$ and
 $\check{\varrho}_\zs{2}=\varrho^{4}_\zs{2}\,\Pi(x^{4})$.


\begin{proposition}\label{Pr.sec:Stc.2-1}
For any
 left continuous functions with finite right limits $f$ and $g$
\begin{equation}\label{sec:Stc.6-01-3}
\E\, V_\zs{t}(f) V_\zs{t}(g)=
\int^{t}_\zs{0}\,
e^{2a(t-s)}
\,
H_\zs{s}(f,g)
\,\d s
\end{equation}
where
$
H_\zs{t}(f,g)=g(t)A_\zs{t}(f)+
f(t)A_\zs{t}(g)
+\sigma_\zs{Q}^2\check{\tau}_\zs{t}(f,g)
+\check{\varrho}_\zs{2}\,
f(t)g(t)$.
\end{proposition}
\proof Applying again \eqref{sec:Stc.2} with $g=1$
yields
\begin{equation}\label{sec:Stc.6-02}
\d V_\zs{t}(f)
=a
V_\zs{t}(f)\d t
+
a
\,f(t)\,\wt{I}_\zs{t}(1)
\d t
+
\d\,L_\zs{t}(f)
\,,
\end{equation}
where
$
L_\zs{t}(f)
=
\int^{t}_\zs{0}\,
\check{I}_\zs{s-}(f)
\d u_\zs{s}
+
\varrho^{2}_\zs{2}
\int^{t}_\zs{0}\,f(s)
\d m_\zs{s}$ and
$\check{I}_\zs{s}(f)=f(s)\xi_\zs{s}
+I_\zs{s}(f)$.
By the Ito formula we get
\begin{align*}
\d V_\zs{t}(f) V_\zs{t}(g)
&= 2a V_\zs{t}(f) \,V_\zs{t}(g)\d t
+
a\left(
g(t)V_\zs{t}(f)
+
f(t)
V_\zs{t}(g)
\right)\,\wt{I}_\zs{t}(1)\d t
\\[2mm]
&
+\d\,[L(f)\,,\,L(g)]_\zs{t} +
V_\zs{t-}(f) \d L_\zs{t}(g)
+
V_\zs{t-}(g) \d L_\zs{t}(f)\,.
\end{align*}

Now from Lemma \ref{Lem.A.0-1} we obtain that
\begin{align}\nonumber
\d \E\,V_\zs{t}(f) V_\zs{t}(g)
&= 2a \E\,V_\zs{t}(f) \,V_\zs{t}(g)\d t
+
\left(
g(t)A_\zs{t}(f)
+
f(t)
A_\zs{t}(g)
\right)\,
\d t
\\[2mm]\label{sec:Stc.6-02-0}
&
+\d\,\E\,[L(f)\,,\,L(g)]_\zs{t}\,,
\end{align}
where $A_\zs{t}(f)=a\E\,V_\zs{t}(f)\,\wt{I}_\zs{t}(1)=a\E\,V_\zs{t}(f)\,V_\zs{t}(1)$.
Note that $\E\,\check{I}_\zs{s}(f)\check{I}_\zs{s}(g)=\sigma_\zs{Q}\check{\tau}_\zs{s}(f,g)$
 and
\begin{align*}
\E\,[L(f)\,,\,L(g)]_\zs{t}
&=\varrho^{2}_\zs{1}\,\int^{t}_\zs{0}\,
\E\,\check{I}_\zs{s}(f)\,
\check{I}_\zs{s}(g)
\,
\d s
+
\E\,\sum_\zs{0\le s\le t}\,\Delta\,L_\zs{s}(f)
\Delta\,L_\zs{s}(g)\\[2mm]
&=\sigma_\zs{Q}^2\,\int^{t}_\zs{0}\,\check{\tau}_\zs{s}(f,g)\,\d s
+
\check{\varrho}_\zs{2}
\,\int^{t}_\zs{0}\,f(s)g(s)\d s
\,.
\end{align*}
To find the function $A_\zs{t}(f)$ we put $g=1$ in \eqref{sec:Stc.6-02-0}.
Taking into account that $A_\zs{t}(1)=\wt{I}^{2}_\zs{t}=\wt{\xi}^{2}_\zs{t}$
we get
$$
\E\,V_\zs{t}(f)V_\zs{t}(1)
=\int^{t}_\zs{0}\,e^{3a(t-s)}\,\left(a\,f(s)\,\E\,\wt{\xi}^{2}_\zs{s} +
\sigma_\zs{Q}^2\,\check{\tau}_\zs{s}(f,1)
+
\check{\varrho}_\zs{2}\,
f(s)
\right)\d s
\,.
$$
Using here that
\begin{equation}\label{sec:Stc.7-06-1-0}
a\tau_\zs{t}(1,1)=(e^{2at}-1)/2
\quad\mbox{and}\quad
a\,\tau_\zs{t}(1,f)=\check{\varepsilon}_\zs{t}(f)
\,,
\end{equation}
we obtain the representation \eqref{sec:Stc.6-01-00-1}. Hence Proposition \ref{Pr.sec:Stc.2-1}.
\endproof

\begin{proposition}\label{Pr.sec:Stc.2}
For any
 left continuous function $f$ with finite right limits
\begin{equation}\label{sec:Stc.7-06-1}
\E\, \wt{I}_\zs{t}(f)\,\wt{I}_\zs{t}(1)
 =
 \int^{t}_\zs{0}\,
e^{2a(t-s)}
\,
\wt{\varkappa}_\zs{s}(f)
\,
\d s\,,
\end{equation}
where
$
\wt{\varkappa}_\zs{s}(f)
=
2
f(s)\,
A_\zs{s}(f)
\,+
4\sigma_\zs{Q}^2\,
f(s)
\tau_\zs{s}(f,1)
+
\check{\varrho}_\zs{2}\,f^{2}(s)
$.
\end{proposition}
\proof
 Using the Ito formula and Lemma \ref{Lem.A.0-1}
we obtain that for any bounded nonrandom functions $f$ and $g$
\begin{align}\nonumber
\d \E\, \wt{I}_\zs{t}(f)\,V_\zs{t}(g)
 &=a \E\, \wt{I}_\zs{t}(f)\,V_\zs{t}(g)\,\d t
+2af(t)\E\,V_\zs{t}(f)\,V_\zs{t}(g)\,\d t\\[2mm]\label{sec:Stc.7-06-2}
&+
a\,g(t)\,\E\,\wt{I}_\zs{t}(f)\,\wt{I}_\zs{t}(1)\d t
+
\d \E\,[\wt{M}(f)\,,\,L(g)]_\zs{t}
\,.
\end{align}
Putting here $g=1$ and taking into account that $V_\zs{t}(1)=\wt{I}_\zs{t}(1)$,
we obtain that
\begin{align*}
\d \E\, \wt{I}_\zs{t}(f)\,V_\zs{t}(1)
 &=2a \E\, \wt{I}_\zs{t}(f)\,V_\zs{t}(1)\,\d t
+2af(t)\E\,V_\zs{t}(f)\,V_\zs{t}(1)\,\d t
\\[2mm]&
+
\d \E\,[\wt{M}(f)\,,\,L(1)]_\zs{t}
\,.
\end{align*}
By the direct calculation we find
$$
\E\,[\wt{M}(f)\,,\,L(1)]_\zs{t}=\int^{t}_\zs{0}\,\check{A}_\zs{s}(f)\d s\,.
$$
So, we get \eqref{sec:Stc.7-06-1} and this proposition.
\endproof

\noindent Further we need the following correlation measures
for two integrated $[0,+\infty)\to \bbr$ functions $f$ and $g$
\begin{equation}\label{sec:Ou.10}
\varpi_\zs{n}(f,g)
=\max_\zs{0\le v+t\le n}\,
\left(
\left|\int^{t}_\zs{0}f(u+v)g(u)\d u\right|
+
\left|\int^{t}_\zs{0} g(u+v) f(u)\d u\right|
\right)
\end{equation}
For any bounded $[0,\infty)\to\bbr$ function $f$ we introduce the following uniform norm
$$
\|f\|_\zs{*,n}=\sup_\zs{0\le t\le n}|f(t)|\,.
$$

\begin{proposition}\label{Pr.sec:Stc.3}
Let  $f$ and $g$ be two
 left continuous bounded by $\phi_\zs{*}$ functions with finite right limits, i.e. $\Vert f\Vert_\zs{*,n}\le \phi_\zs{*}$
 and $\Vert g\Vert_\zs{*,n}\le \phi_\zs{*}$. Then for any $0\le t\le n$
\begin{equation}\label{sec:Stc.7-06-1+1}
\left\vert a\E\,\wt{I}_\zs{t}(f)\,V_\zs{t}(g)
\right\vert\le u^{*}_\zs{1}\varpi_\zs{t}(1,g)
+u^{*}_\zs{2}\varpi_\zs{t}(f,g)
+u^{*}_\zs{3}\,,
\end{equation}
where $u^{*}_\zs{1}=4\phi^{2}_\zs{*}(a_\zs{max})\check{\varrho}_\zs{2}+3\sigma_\zs{Q}^2$,
$u^{*}_\zs{2}=44\phi_\zs{*}\sigma_\zs{Q}^2$ and $u^{*}_\zs{3}=3\phi^{3}_\zs{*}\check{\varrho}_\zs{2}$.
\end{proposition}
\proof
First, note that from Ito formula
we find
\begin{align}\nonumber
a\E\,\wt{I}_\zs{t}(f)\,V_\zs{t}(g)&=
a^{2}\int^{t}_\zs{0}\,e^{a(t-s)}\,g(s)\,\left(\E\,\wt{I}_\zs{s}(f)\,\wt{I}_\zs{s}(1) \right)\d s
\\[2mm] \nonumber
&+2a^{2}
\int^{t}_\zs{0}\,e^{a(t-s)}\,f(s)\,\left(\E\,V_\zs{s}(g)\,V_\zs{s}(f) \right)\d s
\\[2mm]\label{sec:Stc.7-06-1+2}
&+a
\int^{t}_\zs{0}\,e^{a(t-s)}\,
\d \E\,[\wt{M}(f)\,,\,L(g)]_\zs{s}
\,.
\end{align}
Using here
 Lemma \ref{Le.sec:A.7-06-2}.
and
Lemma \ref{Le.sec:A.7-06-3}
we can obtain that
\begin{equation}
\label{sec:Stc.11-06-1}
\vert
a \E\,V_\zs{t}(g)\,V_\zs{t}(f)
\vert
\le 15 \sigma_\zs{Q}^2\,\varpi_\zs{t}(f,g)
+\check{\varrho}_\zs{2}\,\Vert f\Vert_\zs{*,t}
\Vert g\Vert_\zs{*,t}
\,.
\end{equation}

One can check directly that
\begin{align*}
\E\,[\wt{M}(f)\,,\,L(g)]_\zs{t}&=2\sigma_\zs{Q}\,
\int^{t}_\zs{0}\,g(s)f(s)\,
\left(
\E\,I_\zs{s}(f)\,I_\zs{s}(1)
\right)
\d s\\[2mm]
&
+
2\sigma_\zs{Q}\,
\int^{t}_\zs{0}\,f(s)\,
\left(
\E\,I_\zs{s}(f)\,I_\zs{s}(g)
\right)
\d s
+\check{\varrho}_\zs{2}
\int^{t}_\zs{0}\,f^{2}(s)\,g(s)
\d s\,.
\end{align*}
From \eqref{sec:Stc.1} we find that
\begin{align*}
\E\,[\wt{M}(f)\,,\,L(g)]_\zs{s}&=
2\sigma_\zs{Q}^2\,
\int^{t}_\zs{0}\,
g(s)f(s)\,\tau_\zs{s}(f,1)\,
\d s
\\[2mm]
&
+
2\sigma_\zs{Q}^2\,
\int^{t}_\zs{0}\,
f(s)\,\tau_\zs{s}(f,g)\,
\d s
+\check{\varrho}_\zs{2}
\int^{t}_\zs{0}\,f^{2}(s)\,g(s)
\d s\,.
\end{align*}
Using the last equality in \eqref{sec:Stc.7-06-1-0}
we obtain that
\begin{align*}
a\int^{t}_\zs{0}\,e^{a(t-s)}\,
&\d \E\,[\wt{M}(f)\,,\,L(g)]_\zs{s}
=2\sigma_\zs{Q}^2\,\int^{t}_\zs{0}\,e^{a(t-s)}\,
g(s)f(s)\,\check{\varepsilon}_\zs{s}(f)\,
\d s\\[2mm]
&+
2\sigma_\zs{Q}^2\,a\int^{t}_\zs{0}\,e^{a(t-s)}\,
f(s)\,\tau_\zs{s}(f,g)\,
\d s
+
a
\check{\varrho}_\zs{2}
\int^{t}_\zs{0}\,e^{a(t-s)}\,
f^{2}(s)\,g(s)
\d s\,.
\end{align*}
Note now that
$$
\check{\varepsilon}^{\prime}_\zs{t}(f)=a\check{\varepsilon}_\zs{t}(f)
+af(t)(1+e^{2at})/2\,,
$$
i.e.
$\Vert\check{\varepsilon}^{\prime}(f) \Vert_\zs{*,t} \le 2\vert a\vert \Vert f\Vert_\zs{*,t}$. Therefore, in view of  Lemma \ref{Le.sec:A.7-06-1}
we get
$$
\left\vert
\int^{t}_\zs{0}\,e^{a(t-s)}\,
g(s)f(s)\,\check{\varepsilon}_\zs{s}(f)\d s
\right\vert
\le 4\varpi_\zs{t}(f,g) \Vert f\Vert_\zs{*,t}\,.
$$
Moreover, by integrating by parts we can obtain directly that
$$
\left\vert
\int^{t}_\zs{0}\,
g(s)\,\check{\varepsilon}_\zs{s}(f)\d s
\right\vert
\le \varpi_\zs{t}(f,g)\,,
$$
and, therefore,
\begin{equation}
\label{sec:Stc.7-00+1}
\vert \tau_\zs{t}(f,g)\vert\le 3\varpi_\zs{t}(f,g)\,.
\end{equation}
So, the last term in \eqref{sec:Stc.7-06-1+2}
can be estimated as
$$
\left\vert
\,a\int^{t}_\zs{0}\,e^{a(t-s)}\,
\d \E\,[\wt{M}(f)\,,\,L(g)]_\zs{s}
\right\vert
\le 14\,\sigma_\zs{Q}^2
\varpi_\zs{t}(f,g)\,\Vert f\Vert_\zs{*,t}
+\check{\varrho}_\zs{2}\Vert f\Vert^{2}_\zs{*,t} \Vert g\Vert_\zs{*,t}
\,.
$$
Using  Lemma \ref{Le.sec:A.7-06-2-01} in \eqref{sec:Stc.7-06-1+2}
we come to the bound \eqref{sec:Stc.7-06-1+1}. Hence Proposition \ref{Pr.sec:Stc.3}.
\endproof

\begin{proposition}\label{Pr.sec:Stc.4}
Let  $f$ and $g$ be two
 left continuous bounded by $\phi_\zs{*}$ functions with finite right limits, i.e. $\Vert f\Vert_\zs{*,n}\le \phi_\zs{*}$
 and $\Vert g\Vert_\zs{*,n}\le \phi_\zs{*}$. Then for any $t>0$
\begin{equation}\label{sec:Stc.8+1}
\left\vert \E\,[\wt{M}(f)\,,\, \wt{M}(g)]_\zs{t}\,
\right\vert\le
\left(
12\sigma_\zs{Q}^2\,\phi^{2}_\zs{*}\varpi_\zs{t}(f,g)
+\phi^{4}_\zs{*}\check{\varrho}_\zs{2}
\right)
t\,.
\end{equation}
\end{proposition}
\proof
First of all note that from \eqref{Pr.sec:Stc.1}
we obtain that
\begin{align}\nonumber
\E\,[\wt{M}(f)\,,\, \wt{M}(g)]_\zs{t}&=4\sigma_\zs{Q}^2\,\int^{t}_\zs{0}\,f(s)g(s)\tau_\zs{s}(f,g)\,
\d s
\\[2mm] \label{sec:Stc.8+1-1}
&
+
\check{\varrho}_\zs{2}\,\int^{t}_\zs{0}\,f^{2}(s)g^{2}(s)\,\d s\,.
\end{align}
Using here the bound \eqref{sec:Stc.7-00+1}
we obtain  \eqref{sec:Stc.8+1}. Hence Proposition \ref{Pr.sec:Stc.4}.
\endproof

\begin{corollary}\label{Co.sec:Stc.3-1}
Let  $f$ and $g$ be two
 left continuous bounded by $\phi_\zs{*}$ functions with finite right limits, i.e. $\Vert f\Vert_\zs{*,n}\le \phi_\zs{*}$
 and $\Vert g\Vert_\zs{*,n}\le \phi_\zs{*}$. Then for any $t>0$
\begin{equation}\label{sec:Stc.15-06-1+1}
\left\vert\E\,\wt{I}_\zs{t}(f)\wt{I}_\zs{t}(g)\right\vert\,\le\,\left(v^{*}_\zs{1}(\varpi_\zs{t}(1,f)+\varpi_\zs{t}(1,g))+v^{*}_\zs{2}\varpi_\zs{t}(f,g)+v^{*}_\zs{3}\right)\,t\,,
\end{equation}
where $v^{*}_\zs{1}=8\phi^{3}_\zs{*}a_\zs{max}\check{\varrho}_\zs{2}+6\sigma_\zs{Q}^2$,
$v^{*}_\zs{2}=100\phi^2_\zs{*}\sigma_\zs{Q}^2$ and $v^{*}_\zs{3}=13\phi^{4}_\zs{*}\check{\varrho}_\zs{2}$.
\end{corollary}
\proof
From \eqref{sec:Stc.6-01} by the Ito formula one finds for $t\ge 0$
\begin{align}\nonumber
\E\,\wt{I}_\zs{t}(f)\wt{I}_\zs{t}(g)&=\E\,[\wt{M}(f)\,,\, \wt{M}(g)]_\zs{t}
\\[2mm]\label{sec:Stc.4-1-20}
&+
2a\int^{t}_\zs{0}\,\left(f(s)\E\,\wt{I}_\zs{s}(g)\,V_\zs{s}(f)+g(s)\E\,\wt{I}_\zs{s}(f)\,V_\zs{s}(g)\right)\d s\,.
\end{align}
 Using here Proposition \ref{Pr.sec:Stc.3} and Proposition \ref{Pr.sec:Stc.4} we come to desire result.

\endproof

\noindent
Now we set
\begin{equation}\label{sec:Stc.9-00+1}
\overline{I}_\zs{n}(x)=
\sum^{n}_\zs{j=1}\,x_\zs{j}\,\wt{I}_\zs{n}(\phi_\zs{j})\,.
\end{equation}

For this we show the following proposition.

\begin{proposition}\label{Pr.sec:Stc.5}
Assume that $\phi_\zs{1}\equiv 1$.
 Then for any $n\ge 1$
\begin{equation}\label{sec:Stc.9-00+1}
\E\,\overline{I}^{2}_\zs{n}(x)
\le 4 n^{2}\,\left((2\varpi^{*}_\zs{n}+4\phi_\zs{*}^2)v^{*}_\zs{1}
+(\varpi^{*}_\zs{n}+8\phi_\zs{*}^2)v^{*}_\zs{2}
+(1+2\phi_\zs{*}^2)v^{*}_\zs{3}\right),
\end{equation}
where $\varpi^{*}_\zs{n}=\sup_\zs{|i-j|\ge 2}\,\varpi_\zs{n}(\phi_\zs{i},\phi_\zs{j})$.
\end{proposition}

\proof
We
represent the sum as
$$
\sum^n_\zs{j=1}\,x_\zs{j}\,\wt{I}_\zs{n}(\phi_\zs{j})=J_\zs{1,n}+
J_\zs{2,n}\,,
$$
where $J_\zs{1,n}=x_\zs{1}\wt{I}_\zs{n}(\phi_\zs{1})+x_\zs{2}
\wt{I}_\zs{n}(\phi_\zs{2})$ and $J_\zs{2,n}=\sum^n_\zs{j=
3}\,x_\zs{j}\,\wt{I}_\zs{n}(\phi_\zs{j})$. From here we have
\begin{equation}\label{sec:Pr.3-J12}
\E
\overline{I}^{2}_\zs{n}(x)
\,
\le
2\,
\left(
\E\,J^{2}_\zs{1,n}
+
\E\,J^{2}_\zs{2,n}
\right)\,.
\end{equation}
By applying the Cauchy-Schwarz-Bounyakovskii inequality and noting
that $x^{2}_\zs{1}+x^{2}_\zs{2}\le 1$, one gets
$$
\E\,J^{2}_\zs{1,n}\le
\E\,\wt{I}^{2}_\zs{n}(\phi_\zs{1}) +
\E\,\wt{I}^{2}_\zs{n}(\phi_\zs{2}).
$$
Corollary~\ref{Co.sec:Stc.3-1} implies
$$
\E_\zs{Q,S}\,J^{2}_\zs{1,n}\le 4\phi_\zs{*}^2n^2 \left(
2v^{*}_\zs{1}+v^{*}_\zs{2}+v^{*}_\zs{3}
\right)\,.
$$
Here we use that each
$\varpi_\zs{n}(\phi_\zs{i},\phi_\zs{j})\le 2\phi_\zs{*}^2n$.

Applying Corollary~\ref{Co.sec:Stc.3-1}, one gets
\begin{equation}\label{sec:Pr.4-J2}
\E\,J^{2}_\zs{2,n}=2\sum^n_\zs{i,j=
3}\,x_\zs{i}x_\zs{j}\,\E\wt{I}_\zs{n}(\phi_\zs{i})\wt{I}_\zs{n}(\phi_\zs{j})
\le 2n
\sum^n_\zs{i,j= 3} |x_\zs{i}| |x_\zs{j}|
\wt{\kappa}_\zs{i,j}\,,
\end{equation}
where
$\wt{\kappa}_\zs{i,j}=v^{*}_\zs{1}(\varpi_\zs{n}(1,\phi_i)+\varpi_\zs{n}(1,\phi_j))
+v^{*}_\zs{2}\varpi_\zs{n}(\phi_i,\phi_j)+v^{*}_\zs{3}$.
We can estimate the coefficient $\varpi_\zs{i,j}=\varpi_\zs{n}(\phi_i,\phi_j)$
for any $i\ge 3$ as
$\varpi_\zs{i,j}\le 2\phi_\zs{*}^2n\Chi_\zs{\{|i-j|\le 1\}}+
\varpi^{*}_\zs{n}\Chi_\zs{\{|i-j|\ge 2\}}$. By making use of this estimate in
\eqref{sec:Pr.4-J2} and taking into account that
$$
\sum_\zs{i,j\ge 1} |x_\zs{i}| |x_\zs{j}|\le 1
\quad\mbox{and}\quad
\sum_\zs{i,j\ge 3} \Chi_\zs{\{|i-j|\le 1\}}|x_\zs{i}| |x_\zs{j}|\le 3
\,,
$$
 one gets
$$
\sum_\zs{i,j\ge 3} |x_\zs{i}| |x_\zs{j}|
\wt{\kappa}_\zs{i,j}
\,
\le
n\left(2\varpi^*_\zs{n}v^{*}_\zs{1}+(6\phi_\zs{*}^2+\varpi^*_\zs{n})v^{*}_\zs{2}+v^{*}_\zs{3}\right).
$$
From here and the inequalities \eqref{sec:Pr.3-J12}--\eqref{sec:Pr.4-J2} we come to the desired assertion.
Hence Proposition \ref{Pr.sec:Stc.5}. \endproof

Now we check the conditions $\C_1)$ and $\C_2)$ for Ornstein--Uhlenbeck model.
For this we will use the trigonometric basis \eqref{sec:In.5}.

Note that in this case the proxy variance $\sigma_\zs{Q}> 0$ is defined in \eqref{sec:Ex.01-1}.

\begin{proposition}\label{Pr.sec:L1-10}
 Then for any $Q\in\cQ_\zs{n}$ and any $n\ge 1$
 $$
 \L_\zs{1,n}(Q)\le 2\,\sigma^{2}_\zs{Q}(4 a^{2}+15\vert a\vert+2)
 \,.
 $$
\end{proposition}

\proof
First
we note that
\begin{equation}\label{sec:Pr.1-0}
\E_\zs{Q,S} \xi^2_\zs{j,n}= \sigma_\zs{Q} \left( 1+b_\zs{j,n}
\right)\,,
\end{equation}
where
$b_\zs{j,n}=n^{-1}a\int^{n}_\zs{0}\,e^{av}\,\Upsilon_\zs{j}(v)\d
v$ and
$$
 \Upsilon_\zs{j}(v)=\int^{n-v}_\zs{0}\,\Trg_\zs{j}(t+v)\,\Trg_\zs{j}(t)\,\left(
1+e^{2at}
\right)\d t\,.
$$
If $j=1$, one has
\begin{equation}\label{sec:Pr.1}
|\E_\zs{Q,S} \xi^2_\zs{1,n}
-
\sigma_\zs{Q}
 |
\le 2\sigma_\zs{Q}\,.
\end{equation}
Since for the trigonometric basis \eqref{sec:In.5}
for $j\ge 2$
$$
\Trg_\zs{j}(t+v)\,\Trg_\zs{j}(t)=
\cos(\gamma_\zs{j}v)+(-1)^{j}\cos(\gamma_\zs{j}(2t+v))
$$
where $\gamma_\zs{j}=2\pi [j/2]$, therefore,
$$
\Upsilon_\zs{j}(v)=\cos(\gamma_\zs{j}v) F(v)+
(-1)^{j}\,\Upsilon_\zs{0,j}(v)
\,, \quad
F(v)=\int^{n-v}_\zs{0}
\left(
1+e^{2a t}
\right)\d t
$$
 and
$$
\Upsilon_\zs{0,j}(v)=
\int^{n-v}_\zs{0}
\cos(\gamma_\zs{j}(2t+v))
\left(
1+e^{2a t}
\right)\d t\,.
$$
Integrating by parts one finds
$$
\Upsilon_\zs{0,j}(v)
=-\frac{2+e^{2a(n-v)}}{2\gamma_\zs{j}}\,
\sin(v\gamma_\zs{j})+
\frac{a}{2\gamma^{2}_\zs{j}}\,\Upsilon_\zs{1,j}(v)
$$
where
$$
\Upsilon_\zs{1,j}(v)=
\cos(v\gamma_\zs{j}) (e^{2a(n-v)}-1)
-
2a\int^{n-v}_\zs{0} e^{2at} \cos((2t+v)\gamma_\zs{j})\,\d t\,.
$$
It is obvious that $|\Upsilon_\zs{1,j}(v)|\le 2$. Further we
calculate
\begin{align*}
b_\zs{j,n}&=
\frac{a}{n}\int^{n}_\zs{0}\,e^{av}\,F(v)\,\cos(v\gamma_\zs{j})\d v
+
\frac{a}{n}(-1)^{j}
\int^{n}_\zs{0}\,e^{av}\,\Upsilon_\zs{0,j}(v)\,\d v
\\[2mm]
&
:=a D_\zs{1,j}+
a(-1)^{j} D_\zs{2,j}\,.
\end{align*}
Integrating by parts two times yields
$$
D_\zs{1,j}=\frac{1}{n\gamma^{2}_\zs{j}}
\left(
e^{an}\dot{F}(n)
-
\dot{F}(0)-aF(0)-
\int^{n}_\zs{0}\,e^{av} F_\zs{1}(v)\d v
\right)\,,
$$
where $F_\zs{1}(v)=a^{2}F(v)+2a\dot{F}(v)+\ddot{F}(v)$. Since
$\gamma_\zs{j}\ge j$  for $j\ge 2$, we obtain
$$
|D_\zs{1,j}|\le \frac{1}{j^{2}}
\left(
4|a|
+10
\right)\,.
$$
Similarly, one gets $|D_\zs{2,j}|\le 5/j^{2}$. Substituting these
estimates in \eqref{sec:Pr.1-0} and using the upper bound
\eqref{sec:Pr.1}, we obtain for all $j\ge 1$
\begin{equation}\label{sec:Pr.2}
|\E_\zs{Q,S} \xi^2_\zs{j,n}
-
\sigma_\zs{Q}
 |
\le \sigma_\zs{Q}\,\frac{(4 a^{2}+15\,|a|+2)}{j^{2}}\,.
\end{equation}
Thus we arrive at the inequality
$$
\L_\zs{1,n}(Q)\le 2 \sigma_\zs{Q}(4 a^{2}+15\,|a|+2)\,.
$$

\endproof

Proposition \ref{Pr.sec:L1-10} and \eqref{sec:Ex.01-1} -- \eqref{sec:Ex.01-2} imply that the condition $\C_\zs{1})$ holds.

\noindent

\begin{proposition}\label{Pr.sec:L2-11}

For any  $n\ge 1$ and $Q\in\cQ_\zs{n}$
$$
\L_\zs{2,n}(Q)\le 8 \M_\zs{Q}
\,,
$$
where
$\M_\zs{Q}=48\sqrt{2}\vert a\vert
\check{\varrho}_\zs{2}+918 \sigma_\zs{Q}+65\check{\varrho}_\zs{2}$.
\end{proposition}

\proof
We note that for the trigonometric basis \eqref{sec:In.5}
$\|\Trg_\zs{j}\|_\zs{*,n}\le \sqrt{2}$ and $\varpi^{*}_\zs{n}=2$. Indeed,
for any $i\ge 3$,
$$
\Trg_\zs{i}(v+u)=\varkappa_\zs{1,i}(v)\Trg_\zs{i-1}(u)
+\varkappa_\zs{2,i}(v)\Trg_\zs{i}(u)
+\varkappa_\zs{3,i}(v)\Trg_\zs{i+1}(u),
$$
where $\varkappa_\zs{i,j}(\cdot)$ are bounded functions.
 From here in view of the orthonormality and the periodicity of
the functions $(\Trg_\zs{j})_\zs{j\ge 1}$, it follows that for
$0\le t\le n$ and $|i-j|\ge 2$
\begin{align*}
\left|
\int^{t}_\zs{0}
\Trg_\zs{i}(u+v) \Trg_\zs{j}(u)
\d u
\right|
&=
\left|
\int^{\{t\}}_\zs{0}
\Trg_\zs{i}(u+v) \Trg_\zs{j}(u)
\d u
\right|\\[2mm]
&\le \sqrt{\int^{1}_\zs{0}
\Trg^{2}_\zs{i}(u+v)\d u}
=1\,,
\end{align*}
where $\{t\}$ is the fractional part of $t$. Therefore
$\varpi^{*}_\zs{n}\le 2$ if $|i-j|\ge 2$. Thus,
we have that
$\L_\zs{2,n}(Q)\le 8 \M_\zs{Q}$. Hence Proposition \ref{Pr.sec:L2-11}.
\endproof

Proposition \ref{Pr.sec:L2-11} and \eqref{sec:Ex.01-1} -- \eqref{sec:Ex.01-2} imply that the condition $\C_\zs{2})$ holds.


\begin{proposition}\label{Pr.sec:Stc.2}
Let the noise $(\xi_{t})_{t\geq 0}$\ in equation \eqref{sec:In.1}
describes by non-Gaussian Ornstein--Uhlenbeck process \eqref{sec:Ex.1}.
Assume that the basis function satisfy the conditions $\B_\zs{1}$)
and $\B_\zs{2}$).
Then for all $d\geq d_\zs{0}$ the condition $\D_\zs{2}$) holds with $l^{*}_\zs{n}=\underline{\varrho}_\zs{n}(d-6)/2$.
\end{proposition}

\proof
We have
$$
\xi_t=\varrho_1\xi_t^{(1)}+\varrho_2\xi_t^{(2)},\; 0\leq t\leq n,
$$
where $(\xi_\zs{t}^{(1)})_\zs{t\ge 0}$ and $(\xi_\zs{t}^{(2)})_\zs{t\ge 0}$ are independent
Ornstein--Uhlenbeck processes obey the equations
$$
\d\xi_t^{(1)}=a\xi_t^{(1)}dt+\d w_\zs{t}
\quad\mbox{and}\quad
\d\xi_t^{(2)}=a\xi_t^{(2)}dt+\d z_\zs{t}\,.
$$
Moreover, for any square integrated  functions $f$   we set
\begin{equation}\label{sec:Stc.7}
I_\zs{t}^{(1)}(f)=\int_0^t f(s)d\xi_s^{(1)}
\quad\mbox{and}\quad
I_\zs{t}^{(2)}(f)=\int_0^t f(s)d\xi_s^{(2)}
\,.
\end{equation}
Then the matrix $\cD_\zs{n}$ can be rewritten as
$$
\cD_\zs{n}=\varrho_1^2\,\cD_\zs{1,n}+\varrho_2^2\cD_\zs{2,n},
$$
where the $(i,j)$ element of the matrix $\cD_\zs{l,n}$ is defined as
$\E(I_n^{(l)}(\phi_\zs{i})I_n^{(l)}(\phi_\zs{j})|\cG)$. Using the celebrated inequality
of Lidskii and Wieland (see, for example, in
\cite{MarchallOlkin1979}, G.3.a., p.334
)
 we obtain
\begin{equation}\label{eq1.24}
\tr\cD(\cG)-\lambda_{\max}(\cD(\cG))\ge\varrho_1^{2}(
\tr\,\cG_\zs{1,n}-\lambda_{\max}(\cG_\zs{1,n}))
\quad\mbox{a.s.}
\end{equation}
Now,  using
Proposition \ref{Pr.sec:Stc.1}
with $\varrho_\zs{1}=1$ and
$\varrho_\zs{2}=0$
 we obtain that
\begin{equation}\label{eq2.27}
\tr \,\cG_\zs{n}=
\frac{1}{n}\sum_{j=1}^d\E(I_n^{(1)}(\phi_{j}))^2
=d+\sum_{j=1}^d\,b_\zs{j,n}\,,
\end{equation}
where
$$
b_\zs{j,n}
=\frac{a}{n}\int_\zs{0}^n\phi_{j}(t)\int_\zs{0}^{t}
e^{a(t-s)}\phi_{j}(s)(1+e^{2as})\d s \d t\,.
$$
Therefore, setting $\Phi(t,v)=\sum_{j=1}^d\phi_\zs{j}(t)\phi_\zs{j}(t-v)$,
we get
\begin{align*}
\tr \cG_\zs{1,n}&=d+\frac{a}{n}\int_0^n
e^{av}\left(\int_\zs{v}^{n}\Phi(t,v)
\left(
1+e^{2a(t-v)}
\right)
\d t\right)
\d v\\[2mm]
&\ge
d-
2\vert a\vert
\int_0^n
e^{av}
\,\Phi^{*}_\zs{d}(v)
\d v
\,,
\end{align*}
where the function $\Phi^{*}_\zs{d}(\cdot)$
is defined in the condition $\B_\zs{2}$). Taking into account that this function is $1$ - periodic we conclude that
\begin{align*}
\tr \cG_\zs{1,n}&>d-2\vert a\vert
\sum^{n}_\zs{k=1}\,e^{a(k-1)}
\,
\int_\zs{0}^1
\,\Phi^{*}_\zs{d}(v)
\d v
\\[2mm]
&
>d-
\frac{2a_\zs{max}}{1-e^{-a_\zs{max}}}
\int_\zs{0}^1
\,\Phi^{*}_\zs{d}(v)
\d v
=
d-
\frac{1}{2\check{a}}
\int_\zs{0}^1
\,\Phi^{*}_\zs{d}(v)
\d v
\,,
\end{align*}
where $\check{a}$ is given in the condition $\B_\zs{2}$) which implies immediately
$$
\tr
\cD_\zs{1,n}
>\,
d/2\,.
$$
Now, note that
$$
\lambda_{\max}(\cG_\zs{1,n}) =\sup_{\|z\|_\zs{d}=1}n^{-1}\E(I_n^{(1)}(g_\zs{z}))^2
\,,
$$
where $g_\zs{z}=\sum^{d}_\zs{j=1}\,z_\zs{j}\,\phi_\zs{j}$.
Using again
Proposition \ref{Pr.sec:Stc.1} we find
\begin{align*}
\frac{1}{n}\,
\E(I_n^{(1)}(g_\zs{z}))^2&=
\frac{2 a}{n}\,\int_\zs{0}^{n}g_\zs{z}(t)\,
\overline{g}_\zs{z}(t)
\d t+
\frac{1}{n}\,
\int_\zs{0}^{n}\, g^2_\zs{z}(t)\,
\d t\\[2mm]
&=
\frac{a}{n}\,\int_0^n e^{av}\,G_\zs{z}(v)\,\d v
+\int_\zs{0}^1 g^2_\zs{z}(t)dt.
\end{align*}
where $G_\zs{z}(v)=\int_\zs{v}^{n}\,
g_\zs{z}(t)g_\zs{z}(t-v)(1+e^{2a(t-v)})\d t$.
Note now that for all $a\le 0$
$$
\vert G_\zs{z}(v)\vert
\le 2\int^{n}_\zs{0}\,g^{2}_\zs{z}(t)\,\d t
=2 n
\int^{1}_\zs{0}\,g^{2}_\zs{z}(t)\,\d t
\,.
$$
Moreover, taking into account
$$
\int^{1}_\zs{0}\,g^{2}_\zs{z}(t)\,\d t =\sum^{d}_\zs{j=1}\,z^{2}_\zs{j}=1\,,
$$
we obtain that
$\lambda_{\max}(\cG_\zs{1,n})\le 3$. Hence Proposition \ref{Pr.sec:Stc.2}.
\endproof


\section{Proofs}\label{sec:Prf}

\subsection{Proof of Theorem \ref{Th.sec:Imp.1}}

Consider the quadratic error of the estimate \eqref{sec:Imp.11}
\begin{align*}
\|S^{*}_\zs{\gamma}-S\|^2&=\sum^{n}_\zs{j=1}(\gamma(j)\theta^{*}_\zs{j,n}-\theta_\zs{j})^2=
\sum^{d}_\zs{j=1}(\gamma(j)\theta^{*}_\zs{j,n}-\theta_\zs{j})^2+\sum^{n}_\zs{j=d+1}(\gamma(j)\wh{\theta}_\zs{j,n}-\theta_\zs{j})^2\\
&=\sum^{n}_\zs{j=1}(\gamma(j)\wh{\theta}_\zs{j,n}-\theta_\zs{j})^2
+
\c^{2}_\zs{n}
-
2
\c_\zs{n}
\sum^{d}_\zs{j=1}(\wh{\theta}_\zs{j,n}-\theta_\zs{j})\frac{\wh{\theta}_\zs{j,n}}{\|\wt{\theta}_\zs{n}\|_\zs{d}}
\\
&=\|\wh{S}_\zs{\gamma}-S\|^2+
\c^{2}_\zs{n}
-
2
\c_\zs{n}
\sum^{d}_\zs{j=1}(\wh{\theta}_\zs{j,n}-\theta_\zs{j})
\iota_\zs{j}(\wt{\theta}_\zs{n})
\,,
\end{align*}
where  $\iota_\zs{j}(x)=x_\zs{j}/\|x\|_\zs{d}$
for $x=(x_\zs{j})_\zs{1\le j\le d}\in\bbr^{d}$.
Therefore, we can represent the risk for the improved estimator $S^{*}_\zs{\gamma}$ as
$$
\cR_\zs{Q}(S^{*}_\zs{\gamma},S)=\cR_\zs{Q}(\wh{S}_\zs{\gamma},S)+
\c^{2}_\zs{n}
-
2
\c_\zs{n}
\,
\E_\zs{Q,S}\,\sum^{d}_\zs{j=1}(\wh{\theta}_\zs{j,n}-\theta_\zs{j})\,I_\zs{j,n}
\,,
$$
where $I_\zs{j,n}=\E(\iota_\zs{j}(\wt{\theta}_\zs{n})(\wh{\theta}_\zs{j,n}-\theta_j)|\cG_\zs{n})$.
Now, taking into account that the vector $\wt{\theta}_\zs{n}=(\wh{\theta}_\zs{j,n})_\zs{1\le j\le d}$
is the $\cG_\zs{n}$ conditionally Gaussian vector in $\bbr^{d}$   with mean $\wt{\theta}=(\theta_\zs{j})_\zs{1\le j\le d}$
and covariance matrix $n^{-1}\G_\zs{n}$, we obtain
\begin{equation*}
I_\zs{j,n}
=\int_{\mathbb{R}^d}\,\iota_\zs{j}(x)(x-\theta_j)\p(x|\cG_\zs{n})
\d x\,.
\end{equation*}
Here $\p(x|\cG_\zs{n})$ is the conditional
distribution density of the vector $\wt{\theta}_\zs{n}$, i.e.
\begin{equation*}
\p(x|\cG_\zs{n})=\frac{1}{(2\pi)^{d/2}\sqrt{\det\G_\zs{n}}}
\exp\left(-\frac{(x-\theta)'\,\G^{-1}_\zs{n}(x-\theta)}{2}\right)
\,.
 \end{equation*}
Changing the variables by
$u=\G^{-1/2}_\zs{n}(x-\theta)$,  one finds that
\begin{equation}\label{sec:2.2c}
I_\zs{j,n}=\frac{1}{(2\pi)^{d/2}}\sum_{l=1}^{d}
\g_\zs{j,l}\int_{\mathbb{R}^{d}}\tilde{\iota}_\zs{j,n}(u)u_{l}\exp\left(-\frac{\|u\|^{2}_\zs{d}}
{2}\right) \d u\,,
 \end{equation}
where
$\wt{\iota}_\zs{j,n}(u)=\iota_\zs{j}(\G^{1/2}_\zs{n} u+\theta)$
and $\g_\zs{ij}$ denotes the $(i,j)$-th\ element of  $\G^{1/2}_\zs{n}$.
Furthermore, integrating by parts, the integral $I_\zs{j,n}$ can be
rewritten as
\begin{equation*}
I_\zs{j,n}=\sum_{l=1}^{d}\sum_{k=1}^{d}\E\left(\g_\zs{jl}
\,\g_\zs{kl}\,
\frac{\partial \iota_\zs{j}}{\partial
u_k}(u)|_{u=\wt{\theta}_\zs{n}}|\cG_\zs{n}\right)
\,.
\end{equation*}
Now taking into account that  $z^{\prime}Az\leq\lambda_{max}(A)\|z\|^2$
and the condition $\D_\zs{2}$)
we obtain that
\begin{align*}
\Delta_{Q}(S)&=
\c^{2}_\zs{n}
-
2
\c_\zs{n}
n^{-1}\E_\zs{Q,S}\,\left(\frac{\tr \G_\zs{n}}{\|\wt{\theta}_\zs{n}\|_\zs{d}}-
\frac{\wt{\theta}_\zs{n}^{\prime}\G_\zs{n}\wt{\theta}_\zs{n}}{\|\wt{\theta}_\zs{n}\|^3}\right)
\\[2mm]
&
\le \c^{2}_\zs{n}
-
2
\c_\zs{n}
\,l^{*}_\zs{n} n^{-1}\E_\zs{Q,S}\,\frac{1}{\|\wt{\theta}_\zs{n}\|_\zs{d}}
\,.
\end{align*}
Recall, that the $\prime$ denotes the transposition.
Moreover, in view of the Jensen inequality we can estimate the last expectation from below as
$$
\E_\zs{Q,S}\,(\|\wt{\theta}_\zs{n}\|_\zs{d})^{-1}=\E_\zs{Q,S}\,(\|\wt{\theta}+n^{-1/2}\wt{\xi}_\zs{n}\|_\zs{d})^{-1}
\geq
\,(\|\theta\|_\zs{d}+n^{-1/2}\E_\zs{Q,S}\|\wt{\xi}_\zs{n}\|_\zs{d})^{-1}
\,.
$$
Note now that
 the condition
through the inequality
\eqref{sec:In.3}
we  obtain
$$
\E_\zs{Q,S}\|\wt{\xi}_\zs{n}\|^{2}_\zs{d}\le \varkappa_\zs{Q}\,d\,.
$$
So, for $\Vert S\Vert^{2}\le r^{*}_\zs{n}$
$$
\E_\zs{Q,S}\,\|\wt{\theta}_\zs{n}\|^{-1}\geq
\left(r^{*}_\zs{n}+\sqrt{d\varkappa_\zs{Q}/n}\right)^{-1}
$$
and, therefore,
$$
\Delta_{Q}(S)
\le \c^{2}_\zs{n}
-
2
\c_\zs{n}
\frac{l^{*}_\zs{n}}{\left(r^{*}_\zs{n}+\sqrt{d\varkappa_\zs{*}/n}\right)\,n}
=-\c^{2}_\zs{n}
\,.
$$
Hence
Theorem \ref{Th.sec:Imp.1}.
\endproof

\subsection{Proof of Theorem \ref{sec:Mo.Th.1}}

Substituting \eqref{sec:Mo.4} in \eqref{sec:Mo.1} yields for any $\gamma\in\Gamma$
\begin{align}\nonumber
\Er_\zs{n}(\gamma)\,&=\,J_\zs{n}(\gamma)+
2\,\sum^{n}_\zs{j=1}\,\gamma(j)\left(\theta^*_\zs{j,n}\,\wh{\theta}_\zs{j,n}-\frac{\wh{\sigma}_\zs{n}}{n}
-\theta^*_\zs{j,n}\,\theta_\zs{j}\right)\\[2mm]
\label{sec:Mo.11}
&+\,
\|S\|^2-\rho\wh{P}_\zs{n}(\gamma)\,.
\end{align}
Now we set $L(\gamma)=\sum^{n}_\zs{j=1}\,\gamma(j)$,
$$
B_\zs{1,n}(\gamma)=\sum^{n}_\zs{j=1}\,\gamma(j)(\E_\zs{Q}\xi_\zs{j,n}^2-\sigma_\zs{Q})\,,
\quad
B_\zs{2,n}(\gamma)=\sum^{n}_\zs{j=1}\,\gamma(j)\wt{\xi}_\zs{j,n}\,,
$$
$$
M(\gamma)=\frac{1}{\sqrt{n}}\sum^{n}_\zs{j=1}\,\gamma(j)\theta_\zs{j}\xi_\zs{j,n}
\quad\mbox{and}\quad
B_\zs{3,n}(\gamma)=\frac{1}{\sqrt{n}}\sum^{n}_\zs{j=1}\,\gamma(j)g(j)\wh{\theta}_\zs{j,n}\xi_\zs{j,n}\,.
$$
Taking into account the definition \eqref{sec:Mo.5}, we can rewrite \eqref{sec:Mo.11} as
\begin{equation*}
\Er_\zs{n}(\gamma)\,=\,J_\zs{n}(\gamma)+2\frac{\sigma_\zs{Q}-\wh{\sigma}_\zs{n}}{n}L(\gamma)+
2\,M(\gamma)+\frac{2}{n}B_\zs{1,n}(\gamma)
\end{equation*}
\begin{equation}\label{sec:Mo.12}
+2\sqrt{P_\zs{n}(\gamma)}\frac{B_\zs{2,n}(\overline{\gamma})}{\sqrt{\sigma_\zs{Q}n}}-2B_\zs{3,n}(\gamma)\,+\,
\|S\|^2-\rho\wh{P}_\zs{n}(\gamma)
\end{equation}
with $\overline{\gamma}=\gamma/|\gamma |_\zs{n}$. Let $\gamma_0=(\gamma_0(j))_{1\le n}$ be a fixed sequence in $\Gamma$ and $\gamma^*$ be as in \eqref{sec:Mo.6}.
Substituting $\gamma_0$ and $\gamma^*$ in \eqref{sec:Mo.12}, we consider the difference
\begin{align*}
\Er_\zs{n}(\gamma^*)-\Er_\zs{n}(\gamma_0)&\leq 2\frac{\sigma_\zs{Q}-\wh{\sigma}_\zs{n}}{n}L(x)+2M(x)
+\frac{2}{n}B_\zs{1,n}(x)\\[2mm]
&+2\sqrt{P_\zs{n}(\gamma^*)}\frac{B_\zs{2,n}(\overline{\gamma^*})}{\sqrt{\sigma_\zs{Q}n}}
-2\sqrt{P_\zs{n}(\gamma_0)}\frac{B_\zs{2,n}(\overline{\gamma_0})}{\sqrt{\sigma_\zs{Q}n}}\\[2mm]
&-2B_\zs{3,n}(\gamma^*)+2B_\zs{3,n}(\gamma_0)
-\rho\wh{P}_\zs{n}(\gamma^*)
+\rho\wh{P}_\zs{n}(\gamma_0)\,,
\end{align*}
where $x=\gamma^*-\gamma_0$.
Note that $|L(x)|\leq 2\vert \Gamma\vert_\zs{*}$ and $|B_\zs{1,n}(x)|\le \L_\zs{1,n}(Q)$.
Applying the elementary inequality
\begin{equation}\label{sec:Mo.13}
2|ab|\leq \varepsilon a^2+\varepsilon^{-1} b^2
\end{equation}
with any $\varepsilon>0$, we get
$$
2\sqrt{P_\zs{n}(\gamma)}\frac{B_\zs{2,n}(\overline{\gamma})}{\sqrt{\sigma_Q n}}
\leq \varepsilon P_\zs{n}(\gamma)+\frac{B_\zs{2,n}^2(\overline{\gamma})}{\varepsilon\sigma_Q n}
\leq \varepsilon P_\zs{n}(\gamma)+\frac{B^{*}_\zs{2}}{\varepsilon \sigma n}\,,
$$
 where
 $$
B^{*}_\zs{2}
=
\max_\zs{\gamma\in\Gamma}\,
\left(
B_\zs{2,n}^2(\overline{\gamma})
+
B_\zs{2,n}^2(\overline{\gamma}^{2})
\right)
$$
with $\gamma^2=(\gamma_j^2)_\zs{1\le j\le n}$.
  Note that from  definition the function $\L_\zs{2,n}(Q)$ in the condition $\C_2)$
we obtain that
\begin{equation}\label{sec:Mo.13_Ub}
\E_\zs{Q}\,B^{*}_\zs{2}
\le
\sum_\zs{\gamma\in\Gamma}\,
\left(
\E_\zs{Q}B_\zs{2,n}^2(\overline{\gamma})
+
\E_\zs{Q}B_\zs{2,n}^2(\overline{\gamma}^{2})
\right)
\le 2\nu \L_\zs{2,n}(Q)\,.
\end{equation}

Moreover, by the same method we estimate the term $B_\zs{3,n}$.
Note that
\begin{equation}\label{sec:Mo.13_Ub++c-n}
\sum^{n}_\zs{j=1}\,g^{2}_\zs{\gamma}(j)\,
\wh{\theta}^{2}_\zs{j}
=\c^{2}_\zs{n}
\le \frac{\c^{*}_\zs{n}}{n}
\,,
\end{equation}
where $\c^{*}_\zs{n}=n\max_{\gamma\in \Gamma}\c_n^2$.
Therefore,
 through the  Cauchy--Schwarz inequality we can estimate the term $B_\zs{3,n}(\gamma)$ as
$$
|B_\zs{3,n}(\gamma)|\le
\frac{|\gamma |_\zs{n}}{\sqrt{n}}\c_\zs{n}
\left(\sum_\zs{j=1}^{n}
\overline{\gamma}^{2}(j)
\,
\xi^{2}_\zs{j}
\right)^{1/2}
=
\frac{|\gamma |_\zs{n}}{\sqrt{n}}\c_\zs{n}
\left(
\sigma_Q+ B_\zs{2,n}(\overline{\gamma}^{2})
\right)^{1/2}
\,.
$$
So, applying the elementary inequality \eqref{sec:Mo.13} with some arbitrary $\varepsilon>0$, we have
$$
2|B_\zs{3,n}(\gamma)|
\leq \varepsilon P_\zs{n}(\gamma)+\frac{\c^{*}_\zs{n}}{\varepsilon\sigma_Q n}
(\sigma_Q+B^{*}_\zs{2})\,.
$$

Using the bounds above, one has
\begin{multline*}
\Er_\zs{n}(\gamma^*)\leq\Er_\zs{n}(\gamma_0)+\frac{4\vert\Gamma\vert_\zs{n} |\wh{\sigma}_\zs{n}-\sigma_Q|}{n}
+2M(x)+\frac{2}{n}\L_\zs{1,n}(Q)
\\[2mm]
+\frac{2}{\varepsilon}\,\frac{\c^{*}}{n\sigma_Q}(\sigma_Q+B^{*}_\zs{2})
+\frac{2}{\varepsilon}\,
\frac{B^{*}_\zs{2}}{n\sigma_Q}
\\[2mm]
+2 \varepsilon P_\zs{n}(\gamma^*)
+2\varepsilon P_\zs{n}(\gamma_0)
-\rho\wh{P}_\zs{n}(\gamma^*)+\rho\wh{P}_\zs{n}(\gamma_0)\,.
\end{multline*}
The setting $\varepsilon=\rho/4$
and the estimating where this is possible $\rho$ by $1$
 in this inequality
imply
\begin{multline*}
\Er_\zs{n}(\gamma^*)\leq\Er_\zs{n}(\gamma_0)+
 \frac{5\vert\Gamma\vert_\zs{n}|\wh{\sigma}_\zs{n}-\sigma_Q|}{n}
+2M(x)+\frac{2}{n}\L_\zs{1,n}(Q)
\\[2mm]
+\frac{16 (\c^{*}_\zs{n}+1)(\sigma_Q+B^{*}_\zs{2})}{\rho n\sigma_Q}
-\frac{\rho}{2}\wh{P}_\zs{n}(\gamma^*)+\frac{\rho}{2} P_\zs{n}(\gamma_0)
+\rho\wh{P}_\zs{n}(\gamma_0)\,.
\end{multline*}
Moreover, taking into account here that
$$
\vert
\wh{P}_\zs{n}(\gamma_0)
-
P_\zs{n}(\gamma_0)
\vert
\le
\frac{\vert\Gamma\vert_\zs{n}|\wh{\sigma}_\zs{n}-\sigma_Q|}{n}
$$
and that $\rho<1/2$,
we obtain that
\begin{multline}\label{sec:Mo.14}
\Er_\zs{n}(\gamma^*)\leq\Er_\zs{n}(\gamma_0)+
 \frac{6\vert\Gamma\vert_\zs{n}|\wh{\sigma}_\zs{n}-\sigma_Q|}{n}
+2M(x)+\frac{2}{n}\L_\zs{1,n}(Q)
\\[2mm]
+\frac{16 (\c^{*}_\zs{n}+1)(\sigma_Q+B^{*}_\zs{2})}{\rho n\sigma_Q}
-\frac{\rho}{2} P_\zs{n}(\gamma^*)+\frac{3\rho}{2} P_\zs{n}(\gamma_0)\,.
\end{multline}

Now we examine the third term in the right-hand side of this inequality. Firstly we note that
\begin{equation}
\label{upper_bound_M-+01}
2|M(x)|\leq\varepsilon\|S_\zs{x}\|^2+\frac{Z^*}{n\varepsilon}\,,
\end{equation}
where $S_\zs{x}=\sum^{n}_\zs{j=1}\,x_\zs{j}\theta_\zs{j}\phi_\zs{j}$
and
$$
Z^*=\sup_{x\in\Gamma_1}\frac{nM^2(x)}{\|S_x\|^2}
\,.
$$
We remind that the set $\Gamma_\zs{1}=\Gamma-\gamma_\zs{0}$.
Using Proposition
\ref{Pr.sec:Stc.1}
we can obtain that for any fixed $x=(x_\zs{j})_\zs{1\le j\le n}\in\bbr^{n}$
\begin{equation}
\label{M^2+11-00}
\E\,M^2(x)=\frac{\E\,I^{2}_\zs{n}\left(S_\zs{x}\right)}{n^{2}}
=\frac{\sigma_Q \Vert S_\zs{x}\Vert^{2}}{n}
=\frac{\sigma_Q}{n}\,
\sum^{n}_\zs{j=1}\,x^{2}_\zs{j}\,\theta^{2}_\zs{j}
\end{equation}
and, therefore,
\begin{equation}
\label{up-Z*-00}
\E_\zs{Q}Z^*
\le
\sum_{x\in\Gamma_1}\frac{n M^2(x)}{\|S_x\|^2}
\leq \sigma_Q
\nu
\,.
\end{equation}
Moreover,  the norm $\Vert S^{*}_\zs{\gamma^{*}}-S^{*}_\zs{\gamma_\zs{0}}\Vert$ can be estimated from below as
\begin{align*}
\Vert S^{*}_\zs{\gamma}-S^{*}_\zs{\gamma_\zs{0}}\Vert^{2}
&=
\sum^{n}_\zs{j=1}
(x(j)+\beta(j))^{2}\wh{\theta}^{2}_\zs{j}
\\[2mm]
&\ge
\|\wh{S}_\zs{x}\|^2
+2
\sum^{n}_\zs{j=1}
x(j)\beta(j)\wh{\theta}^{2}_\zs{j}\,,
\end{align*}
where $\beta(j)=\gamma_0(j)g_j(\gamma_0)-\gamma(j)g_j(\gamma)$.
Therefore, in view of \eqref{sec:Imp.4}
\begin{align*}
\|S_\zs{x}\|^2&-\|S^{*}_\zs{\gamma}-S^{*}_\zs{\gamma_\zs{0}}\|^2
\le
\|S_\zs{x}\|^2-\|\wh{S}_\zs{x}\|^2
-2\sum^{n}_\zs{j=1}\,x(j)\beta(j)
\wh{\theta}_\zs{j}^2
\\[2mm]
&
\le
-2M(x^{2})-2\sum^{n}_\zs{j=1}\,x(j)\beta(j)
\wh{\theta}_\zs{j}\theta_j-
\frac{2}{\sqrt{n}}\Upsilon(x)
\,,
\end{align*}
where $\Upsilon(\gamma)=\sum^{n}_\zs{j=1}\,\gamma(j)\beta(j)
\wh{\theta}_\zs{j}\xi_\zs{j}$. Note that the first term in this inequality we can estimate as
$$
2M(x^{2})\le \varepsilon\|S_\zs{x}\|^2+\frac{Z_1^*}{n\varepsilon}
\quad\mbox{and}\quad
Z^*_1=\sup_{x\in\Gamma_1}\frac{n M^2(x^{2})}{\|S_x\|^2}
\,.
$$
Note that, similarly to \eqref{up-Z*-00} we can estimate the last term as
$$
\E_\zs{Q}Z_1^*\leq \sigma_Q \nu\,.
$$
From this it follows that for any $0<\varepsilon<1$
\begin{align}
\nonumber
\|S_\zs{x}\|^2
\le
\frac{1}{1-\varepsilon}
&\left(
\|S^{*}_\zs{\gamma}-S^{*}_\zs{\gamma_\zs{0}}\|^2
+\frac{Z_1^*}{n\varepsilon}\right.
 \\[2mm] \label{upper-bound-000}
&
\left.
-2\sum^{n}_\zs{j=1}\,x(j)\beta(j)
\wh{\theta}_\zs{j}\theta_j-
\frac{2\Upsilon(x)}{\sqrt{n}}
\right)
\,.
\end{align}
Moreover, note now that the property
 \eqref{sec:Mo.13_Ub++c-n} yields
\begin{equation}
\label{theta-whg-upperb-00}
\sum^{n}_\zs{j=1}\,
\beta^{2}(j)
\wh{\theta}^{2}_\zs{j}
\le
2
\sum^{n}_\zs{j=1}\,g^{2}_\zs{\gamma}(j)\,
\wh{\theta}^{2}_\zs{j}
+
2\sum^{n}_\zs{j=1}\,g^{2}_\zs{\gamma_\zs{0}}(j)\,
\wh{\theta}^{2}_\zs{j}
\le
\frac{4\c^{*}}{\varepsilon n}\,.
\end{equation}
Taking into account that $\vert x(j)\vert\le 1$ and using the inequality \eqref{sec:Mo.13}, we get that
 for any $\varepsilon>0$
$$
2\left|\sum^{n}_\zs{j=1}\,x(j)
\beta(j)
\wh{\theta}_\zs{j}\theta_j\right|\leq\varepsilon\|S_\zs{x}\|^2
+\frac{4\c^{*}}{\varepsilon n}
\,.
$$
To estimate the last term in the right hand of  \eqref{upper-bound-000} we use first
the Cauchy -- Schwarz inequality
and then the bound
\eqref{theta-whg-upperb-00}, i.e.
\begin{align*}
\frac{2}{\sqrt{n}}\vert \Upsilon(\gamma)\vert
&\le
\frac{2\vert\gamma\vert_\zs{n}}{\sqrt{n}}\left(\sum^{n}_\zs{j=1}\,\beta^{2}(j)
\wh{\theta}^{2}_\zs{j}\right)^{1/2}
\left(
\sum^{n}_\zs{j=1}\,\bar{\gamma}^{2}(j)\,\xi^{2}_\zs{j} \right)^{1/2}\\[2mm]
&
\le
\varepsilon P_\zs{n}(\gamma)
+
\frac{\c^{*}}{n\varepsilon\sigma_Q}
\sum^{n}_\zs{j=1}\,\bar{\gamma}^{2}(j)\,\xi^{2}_\zs{j}
\le
\varepsilon P_\zs{n}(\gamma)
+
\frac{\c^{*}(\sigma_Q+B^{*}_\zs{2})}{n\varepsilon\sigma_Q}
\,.
\end{align*}
Therefore,
\begin{align*}
\frac{2}{\sqrt{n}}\vert \Upsilon(x)\vert
&\le
\frac{2}{\sqrt{n}}\vert \Upsilon(\gamma^{*})\vert
+\frac{2}{\sqrt{n}}\vert \Upsilon(\gamma_\zs{0})\vert\\[2mm]
&\le
\varepsilon P_\zs{n}(\gamma^{*})
+\varepsilon P_\zs{n}(\gamma_\zs{0})
+
\frac{2\c^{*}(\sigma_Q+B^{*}_\zs{2})}{n\varepsilon\sigma_Q}
\,.
\end{align*}
So, using all these bounds in \eqref{upper-bound-000},
we obtain that
\begin{align*}
\|S_\zs{x}\|^{2}
&\le
\frac{1}{(1-\varepsilon)}\biggl(\frac{Z_1^*}{n\varepsilon}+
\|S_\zs{\gamma^*}^*-S_\zs{\gamma_0}^*\|^2
+\frac{6\c^{*}_\zs{n}(\sigma+B^{*}_\zs{2})}{n\sigma\varepsilon}
\\[2mm]
&+\varepsilon P_\zs{n}(\gamma^*)+\varepsilon P_\zs{n}(\gamma_0)\biggr)\,.
\end{align*}
Using in the inequality \eqref{upper_bound_M-+01} this bound and the estimate
$$
\|S_\zs{\gamma^*}^*-S_\zs{\gamma_0}^*\|^2\leq
2(\Er_\zs{n}(\gamma^*)+\Er_\zs{n}(\gamma_0))\,,
$$
we obtain
\begin{align*}
2|M(x)|&\le
\frac{Z^*+Z_1^*}{n(1-\varepsilon)\varepsilon}
+
\frac{2\varepsilon(\Er_\zs{n}(\gamma^*)+\Er_\zs{n}(\gamma_0))}{(1-\varepsilon)}
\\[2mm]
&
+
\frac{6\c^{*}_\zs{n}(\sigma_Q+B^{*}_\zs{2})}{n\sigma_Q(1-\varepsilon)}
+\frac{\varepsilon^{2}}{1-\varepsilon}\left( P_\zs{n}(\gamma^*)+ P_\zs{n}(\gamma_0)\right)
\,.
\end{align*}
Choosing here $\varepsilon\le \rho/2<1/2$ we obtain that
\begin{align*}
2|M(x)|&\le
\frac{2(Z^*+Z_1^*)}{n\varepsilon}
+
\frac{2\varepsilon(\Er_\zs{n}(\gamma^*)+\Er_\zs{n}(\gamma_0))}{(1-\varepsilon)}
\\[2mm]
&
+
\frac{12\c^{*}_\zs{n}(\sigma_Q+B^{*}_\zs{2})}{n\sigma_Q}
+
\varepsilon\left( P_\zs{n}(\gamma^*)+ P_\zs{n}(\gamma_0)\right)
\,.
\end{align*}
From here and \eqref{sec:Mo.14}, it follows that
\begin{align*}
\Er_\zs{n}(\gamma^*) &\leq\frac{1+\varepsilon}{1-3\varepsilon}\Er_\zs{n}(\gamma_0)
+ \frac{6\vert\Gamma\vert_\zs{n}|\wh{\sigma}_\zs{n}-\sigma_Q|}{n(1-3\varepsilon)}
+\frac{2}{n(1-3\varepsilon)}\L_\zs{1,n}(Q)
\\[2mm]
&
+\frac{28(1+\c^{*}_\zs{n})(B^{*}_\zs{2}+\sigma_Q)}{\rho(1-3\varepsilon)n\sigma_Q}
+\frac{2(Z^*+Z_1^*)}{n(1-3\varepsilon)}
+\frac{2\rho P_\zs{n}(\gamma_0)}{1-3\varepsilon}.
\end{align*}
Choosing here $\varepsilon=\rho/3$
and estimating $(1-\rho)^{-1}$ by $2$ where this is possible,
 we get
\begin{align*}
\Er_\zs{n}(\gamma^*) &\leq\frac{1+\rho/3}{1-\rho}\Er_\zs{n}(\gamma_0)
+ \frac{12\vert\Gamma\vert_\zs{n}|\wh{\sigma}_\zs{n}-\sigma_Q|}{n}
+\frac{4}{n}\L_\zs{1,n}(Q)
\\[2mm]
&
+\frac{56(1+\c^{*}_\zs{n})(B^{*}_\zs{2}+\sigma_Q)}{\rho n\sigma_Q}
+\frac{4(Z^*+Z_1^*)}{n}
+\frac{2\rho P_\zs{n}(\gamma_0)}{1-\rho}\,.
\end{align*}
Taking the expectation and using the upper bound for $P_\zs{n}(\gamma_0)$  in Lemma~\ref{Lem.A.1} with $\varepsilon=\rho$ yields
$$
\cR_\zs{Q}(S^*,S)\leq\frac{1+5\rho}{1-\rho}\cR_\zs{Q}(S^*_\zs{\gamma_0},S)
+\frac{\check{\U}_\zs{Q,n}}{n\rho}
+
 \frac{12\vert\Gamma\vert_\zs{n}\E_\zs{Q}|\wh{\sigma}_\zs{n}-\sigma_Q|}{n}
\,,
$$
where
$
\check{\U}_\zs{Q,n}=4\L_{1,n}(Q)+
56(1+c^{*}_\zs{n})
(2\L_{2,n}(Q)\nu+1)
+2\c^{*}_\zs{n}$.
The inequality holds for each $\gamma_0\in\Lambda$, this  implies  Theorem \ref{sec:Mo.Th.1}. \endproof

\endproof

\bigskip

\subsection{Proof of Theorem \ref{Th.sec: Ae.2}}

Firstly, note, that for any fixed $Q\in \cQ_\zs{n}$
\begin{equation}\label{sec:Lo.1-0}
\sup_\zs{S\in W_\zs{k,\r}}\,\cR^{*}_\zs{n}(\wh{S}_\zs{n},S)
\ge
\,
\sup_\zs{S\in W_\zs{k,\r}}\,
\cR_\zs{Q}(\wh{S}_\zs{n},S)
\,.
\end{equation}
Now for any fixed $0<\ve<1$ we set
\begin{equation}\label{sec:Lo.1}
d=d_\zs{n}=\left[\frac{k+1}{k}v^{1/(2k+1)}_\zs{n}l_\zs{k}(\r_\zs{\ve})\right]
\quad\mbox{and}\quad
\r_\zs{\ve}=(1-\ve)\r\,.
\end{equation}
Next we approximate
the unknown function by a trigonometric series with $d=d_\zs{n}$ terms, i.e.
for any array $z=(z_\zs{j})_\zs{1\le j\le d_\zs{n}}$,
 we set
\begin{equation}\label{sec:Lo.4}
S_\zs{z}(x)=\sum_{j=1}^{d_\zs{n}}\,z_\zs{j}\,\phi_\zs{j}(x)
\,.
\end{equation}

\noindent
To define the bayesian risk we  choose a prior distribution on $\bbr^{d}$
as
\begin{equation}\label{sec:Lo.5}
\kappa=(\kappa_\zs{j})_\zs{1\le j\le d_\zs{n}}
\quad\mbox{and}\quad
\kappa_\zs{j}=s_\zs{j}\,\eta_\zs{j}\,,
\end{equation}
where $\eta_\zs{j}$ are i.i.d. Gaussian $\cN(0,1)$ random variables and
the coefficients
$$
s_\zs{j}=\sqrt{\frac{s^*_\zs{j}}{v_\zs{n}}}
\quad\mbox{and}\quad
s^{*}_\zs{j}\,
=
\left( \frac{d_\zs{n}}{j}
\right)^{k}
-
1
\,.
$$
Furthermore, for any function $f$, we denote by $\p(f)$ it's projection
in $\L_\zs{2}[0,1]$
 onto
 $W_\zs{k,\r}$, i.e.
 $$
\p(f)=\hbox{\rm Pr}_\zs{W_\zs{k,\r}}(f)\,.
$$
Since $W_\zs{k,\r}$ is a convex set, we obtain
$$
\|\wh{S}-S\|^2\ge\|\wh{\p}-S\|^2
\quad\mbox{with}\quad
\wh{\p}=\p(\wh{S})
\,.
$$
Therefore,
$$
\sup_\zs{S\in W_\zs{k,\r}}\,\cR(\wh{S},S)
\ge\,
\int_{\{z\in\bbr^d\,:\,S_\zs{z}\in W_\zs{k,\r}\}}\,
\E_\zs{S_\zs{z}}\|\wh{\p}-S_\zs{z}\|^2\,\mu_{\kappa}(\d z)
\,.
$$
Using the distribution $\mu_\zs{\kappa}$ we introduce
 the following Bayes risk
$$
\wt{\cR}_\zs{Q}(\wh{S})=
\int_\zs{\bbr^d}\cR_\zs{Q}(\wh{S},S_\zs{z})\,
\mu_\zs{\kappa}(\d z)\,.
$$
Taking into account now that $\|\wh{\p}\|^2\le \r$
 we obtain
\begin{equation}\label{sec:Lo.12}
\sup_\zs{S\in W_\zs{k,\r}}\,
\cR_\zs{Q}(\wh{S},S)
\,
\ge\,
\wt{\cR}_\zs{Q}(\wh{\p})
-2\,
\R_\zs{0,n}
\end{equation}
with
$$
\R_\zs{0,n}=
\int_\zs{ \{z\in\bbr^d\,:\,S_\zs{z}\notin W_\zs{k,\r}\}}\,
\,
(\r+\|S_\zs{z}\|^2)\,
\mu_\zs{\kappa}(\d z)
\,.
$$
Therefore, in view of \eqref{sec:Lo.1-0}
\begin{equation}\label{sec:Lo.1+10}
\sup_\zs{S\in W_\zs{k,\r}}\,\cR^{*}_\zs{n}(\wh{S}_\zs{n},S)
\ge
\,
\sup_\zs{Q\in\cQ_\zs{n}}\,
\wt{\cR}_\zs{Q}(\wh{\p})
-2\,
\R_\zs{0,n}
\,.
\end{equation}
In
Lemma \ref{Le.sec:App.3+1}
we studied the last term in this inequality.
Now it is easy to see that
$$
\|\wh{\p}-S_z\|^2 \ge
\sum_{j=1}^{d_\zs{n}}\,
(\wh{z}_\zs{j}-z_\zs{j})^2
\,,
$$
where $\wh{z}_\zs{j}=\int^{1}_\zs{0}\,\wh{\p}(t)\,\phi_\zs{j}(t)\d t.$ So, in view of Lemma~\ref{Le.sec:App.3}
and reminding
that
$v_\zs{n}=n/\varsigma^{*}$
 we obtain
\begin{align*}
\sup_\zs{Q\in\cQ_\zs{n}}
\wt{\cR}_\zs{Q}(\wh{\p})\,
&\ge\,
\sup_\zs{0<\varrho^{2}_\zs{1}\le \varsigma^{*}}
\sum_{j=1}^{d_\zs{n}}\,\frac{1}
{n\varrho^{-2}_\zs{1}+v_\zs{n}\,(s^{*}_\zs{j})^{-1}}
\\[2mm]
&
=
\frac{1}{v_\zs{n}}\,
\sum_{j=1}^{d_\zs{n}}\,\frac{s^{*}_\zs{j}}
{s^{*}_\zs{j}+\,1}
=
\frac{1}{v_\zs{n}}\,
\sum_{j=1}^{d_\zs{n}}\,
\left(
1
-
\frac{j^k}{d^k_\zs{n}}
\right)
\,.
\end{align*}
Therefore, using now the definition \eqref{sec:Lo.1},
Lemma \ref{Le.sec:App.3+1} and the inequality
\eqref{sec:Lo.1+10}
we obtain that
$$
\liminf_\zs{n\to\infty}\inf_\zs{\wh{S}\in\Sigma_\zs{n}}\,v^{\frac{2k}{2k+1}}_\zs{n}\,
\sup_\zs{S\in W_\zs{k,\r}}\,\cR^{*}_\zs{n}(\wh{S}_\zs{n},S)
\ge\,
(1-\ve)^{\frac{1}{2k+1}}\,
l_\zs{k}(\r_\zs{\ve})\,.
$$
Taking here limit as $\ve\to 0$ implies Theorem~\ref{Th.sec: Ae.2}.
\endproof

\bigskip

\subsection{Proof of Theorem \ref{Th.sec: Ae.1}}

This theorem follows from Theorems \ref{Th.sec:2.3} and \ref{Th.sec:Imp.1} and Theorem 3.1 in \cite{KoPe2009b}.

\endproof

\bigskip
\bigskip

{\bf Acknownledgements.}
This work was partially supported  by the research project no. 2.3208.2017/4.6
(the  Ministry of Education and Science of the Russian Federation) and RFBR Grant 16-01-00121.
The work of the second author was  partially supported
 by the Russian Federal Professor program (project no. 1.472.2016/1.4,
 the  Ministry of Education and Science of the Russian Federation).

\newpage

\medskip

\setcounter{section}{0}
\renewcommand{\thesection}{\Alph{section}}

\section{Appendix}\label{sec:A}

\subsection{Property of the trigonometric basis}

\begin{lemma}\label{Le.sec:A.0}
The
 trigonometric basis \eqref{sec:In.5}
satisfies the conditions $\B_\zs{1}$)
and
$\B_\zs{2}$) with $d_\zs{0}$ and $\check{a}$ defined in \eqref{sec:In.5-00}.
\end{lemma}
\proof
First, we set
$$
\Phi_\zs{d}(t,v)=\sum^{d}_\zs{l=1}\,\phi_\zs{l}(t)\phi_\zs{l}(t-v)
\quad\mbox{and}\quad
N=[d/2]\,.
$$
Note now that for any $d\ge 3$
this sum can be represented as
\begin{align*}
\Phi_\zs{d}(t,v)&
=1+\sum_{j=1}^N\bigl(\phi_{2j}(t)\phi_{2j}(t-v)
+\phi_{2j+1}(t)\phi_{2j+1}(t-v)\bigr)\\[2mm]
&-
\phi_\zs{d+1}(t)\phi_\zs{d+1}(t-v)\,\Chi_\zs{\{d=2N\}}
\\[2mm]
&
=1+2\sum_{j=1}^N\cos(2\pi jv)
-
\phi_\zs{d+1}(t)\phi_\zs{d+1}(t-v)\,\Chi_\zs{\{d=2N\}}\,.
\end{align*}
Therefore,
$$
\Phi^{*}_\zs{d}(v)
\le 2+\vert 1+2\sum_{j=1}^N\cos(2\pi jv)\vert
=2+
\left\vert
\frac{\sin(\pi (2N+1)v)}{\sin(\pi v)}
\right\vert
\,.
$$
So, taking into account that $\vert 1+2\sum_{j=1}^N\cos(2\pi jv)\vert\le 2N+1$,
we obtain that
 for any fixed $0<\delta<1/2$
\begin{align*}
\int^{1}_\zs{0}\,\Phi^{*}_\zs{d}(v)\,\d v
&\le
2+2\delta (d+1)+
\int^{1-\delta}_\zs{\delta}\,\frac{1}{\sin(\pi v)}\,\d v\\[2mm]
&=
2+2\delta (d+1)+2
\int^{1/2}_\zs{\delta}\,\frac{1}{\sin(\pi v)}\,\d v
\,.
\end{align*}
Using here that $\sin(\pi v)\ge 2 v$ for any $0<v<1/2$ we obtain that
$$
\int^{1}_\zs{0}\,\Phi^{*}_\zs{d}(v)\,\d v\le 4+2\delta d
-\ln(2\delta)\,.
$$
Minimizing this upper bound, we obtain that (for $\delta=1/(2d)$)
$$
\int^{1}_\zs{0}\,\Phi^{*}_\zs{d}(v)\,\d v\le 5+\ln d\,.
$$
Hence Lemma \ref{Le.sec:A.0}. \endproof

\subsection{Technical lemmas}

In the following lemma we need   the well-known Novikov  inequalities \cite{Novikov1975}  for purely discontinuous martingales.
Namely, for any $p\ge 2$ and for any $t>0$
\begin{equation}
\label{sec:A.1}
\E\sup_\zs{u\le t}|J_\zs{u}(h)|^p\le \check{C}_\zs{p}
\E\,\big (|h|^{2}*\nu_\zs{t}\big)^{p/2}+\E\, |h|^p*\nu_\zs{t}\,,
\end{equation}
where $J_\zs{t}(h)=g*(\vert\Gamma\vert_\zs{*}-\nu)_t$.

\begin{lemma}\label{Lem.A.0-1}
Assume that $\Pi(x^{8})<\infty$. Then
for any non random left continuous $\bbr_\zs{+}\to \bbr$ functions $f$ and $g$
having the finite right limits and for any $t>0$
$$
\E\,\int^{t}_\zs{0}\,V_\zs{s-}(f)\d L_\zs{s}(g)=0\,.
$$
\end{lemma}
\proof
By the definition we have
\begin{align*}
\int^{t}_\zs{0}\,V_\zs{s-}(f)\d L_\zs{s}(g)&=
\int^{t}_\zs{0}\,V_\zs{s-}(f)
(I_\zs{s-}(g)+g(s)I_\zs{s-}(1))\,\d u_\zs{s}\\[2mm]
&
+\varrho^{2}_\zs{2}\,\int^{t}_\zs{0}\,V_\zs{s-}(f)\,g(s)\d m_\zs{s}
\,.
\end{align*}
Note that in this case $<m>_\zs{t}=\Pi(x^{4})t$. So, to prove this lemma we need to check that
$$
\int^{t}_\zs{0}\,
\E\,V^{2}_\zs{s}(f)\,I^{2}_\zs{s}(g)\d s<\infty
\quad\mbox{and}\quad
\int^{t}_\zs{0}\,
\E\,V^{2}_\zs{s}(f)\,\d s<\infty
\,.
$$
It is clear that to obtain these properties it suffices to show that  for any bounded on the interval $[0,t]$
nonrandom
unction $f$
\begin{equation}
\label{sec:A.0-1}
\sup_\zs{0\le s\le t}\,\E\,I^{8}_\zs{s}(f)\,<\infty
\,.
\end{equation}
Taking into account the definitions in
\eqref{sec:Stc.7} we get
$$
I_\zs{s}(f)=I^{(1)}_\zs{s}(f)
+
I^{(2)}_\zs{s}(f)
\,.
$$
Since the first term $I^{(1)}_\zs{s}(f)$ is the Gaussian random variable
we need to show the inequality
\eqref{sec:A.0-1}
 only for $I^{(2)}_\zs{s}(f)$.
 But it follows immediately from the Novikov inequality
 \eqref{sec:A.1}. Hence lemma \ref{Lem.A.0-1}.

\endproof

\bigskip

\begin{lemma}\label{Le.sec:A.7-06-1}
Let $\upsilon$ be a continuously differentiable $\bbr\to\bbr$ function. Then, for any
$t>0$, $\alpha>0$ and for any  integrated $\bbr\to\bbr$ function $h$,
$$
\left\vert
\int^{t}_\zs{0}\,e^{-\alpha(t-s)} \,h(s) \upsilon(s)\,\d s
\right\vert\,
\le \,
\|\overline{h}\|_\zs{*,t}
\left(
2 \|\upsilon\|_\zs{*,t}
+\frac{\|\dot{\upsilon}\|_\zs{*,t}}{\alpha}\,
\right)\,,
$$
where $\overline{h}_\zs{s}=\int^{t}_\zs{0}\,h(u)\d u$.
\end{lemma}
\proof This  Lemma \ref{Le.sec:A.7-06-1} follows immediately from the integrating par parts.
\endproof

\bigskip

\begin{lemma}\label{Le.sec:A.7-06-2}
For any measurable bounded  $[0,+\infty)\to\bbr$ functions $f$ and $g$,
for any $-\infty< a\le 0 $ and for any $t>0$
$$
\left|
a
\int^{t}_\zs{0}
e^{2a(t-s)}\,g(s)
A_\zs{s}(f)
\,\d s
\right|\,\le
3\,\varrho^{2}_\zs{*}\,
\varpi^{*}_\zs{t}(f,g)
+
\check{\varrho}_\zs{2}
\|f\|_\zs{*,t}
\|g\|_\zs{*,t}
\,,
$$
where
$\check{\varrho}_\zs{2}=\varrho^{2}_\zs{2}\,\Pi(x^{4})$.
\end{lemma}
\proof
First note that
\begin{equation}
\label{sec:Stc.6+00}
A_\zs{t}=J_\zs{t}(f)+\varrho^{2}_\zs{*}\,\check{J}_\zs{t}(f)
\end{equation}
where
$$
J_\zs{t}(f)=\int^{t}_\zs{0}\,e^{3a(t-u)}\,f(u)\upsilon(u)\d u
\quad\mbox{and}\quad
\check{J}_\zs{t}(f)=2\int^{t}_\zs{0}\,e^{3a(t-u)}\,\check{\varepsilon}_\zs{u}(f)\d u
\,.
$$
To study these integrals we need to calculate $\E\wt{\xi}^{2}_\zs{t}$.
To this end
through the equation
\eqref{sec:Stc.6-01}
we can represent this expectation in the following integral form
$$
\E\,\wt{\xi}^{2}_\zs{t}
=
\int^{t}_\zs{0} e^{4a(t-s)}\,
\d \E\,[\wt{M}(1),\wt{M}(1)]_\zs{s}\,.
$$
Moreover, using the definition of $M_\zs{s}(1,1)$
in
\eqref{sec:Stc.2}
we obtain that
$$
\E\,[\wt{M}(1),\wt{M}(1)]_\zs{t}=4\varrho_\zs{*}\,
\int^{t}_\zs{0}\,\E\,\xi^{2}_\zs{s}\d s
+\check{\varrho}_\zs{2}t
\,,
$$
where $\check{\varrho}_\zs{2}=\varrho^{4}_\zs{2}\,\Pi(x^{4})$.
Therefore,
$$
\E\,\wt{\xi}^{2}_\zs{t}
=
\int^{t}_\zs{0} e^{4a(t-s)}\,
\left(
4\varrho_\zs{*}\,\E\,\xi^{2}_\zs{s}
+\check{\varrho}_\zs{2}
\right)
\d s\,.
$$
Using here
\eqref{sec:Stc.3}
we obtain that for $a<0$
\begin{equation}\label{sec:A.4}
\E\,\wt{\xi}^{2}_\zs{t}
=
e^{4at}\frac{2\varrho^{2}_\zs{*}+a\check{\varrho}_\zs{2}}{4a^{2}}
-e^{2at}
\frac{\varrho^{2}_\zs{*}}{a^{2}}
+
\frac{2\varrho^{2}_\zs{*}-a\check{\varrho}_\zs{2}}{4a^{2}}
\,.
\end{equation}
Note now that  the function $\upsilon(\cdot)$ defined in
\eqref{sec:Stc.6-01-00-1}
can be represented as
\begin{equation}\label{sec:A.5}
\upsilon(s)=
a\upsilon_\zs{1}(s)
+
\upsilon_\zs{2}(s)
\end{equation}
with
$$
\upsilon_\zs{1}(s)=
\frac{\check{\varrho}_\zs{2}}{4}
\left(
e^{4at}+3
\right)
\quad\mbox{and}\quad
\upsilon_\zs{2}(s)=\frac{\varrho^{2}_\zs{*}}{2}
\left(
e^{4at}
-
1
\right)
\,.
$$
It is clear that
\begin{equation}\label{sec:A.8-06-1}
\Vert \upsilon_\zs{1}\Vert_\zs{*,n}\le \check{\varrho}_\zs{2}
\quad\mbox{and}\quad
2 \|\upsilon_\zs{2}\|_\zs{*,n}
+\frac{\|\dot{\upsilon}_\zs{2}\|_\zs{*,n}}{2\vert a\vert}
\le 2\varrho_\zs{*}\,.
\end{equation}
Now we  have
$$
J_\zs{t}(f)=J_\zs{1,t}(f)
+
J_\zs{2,t}(f)\,,
$$
where $J_\zs{1,t}(f)=a\int^{t}_\zs{0}\,e^{3a(s-u)}\,f(u)\upsilon_\zs{1}(u)\d u$
and $J_\zs{2,t}(f)=\int^{t}_\zs{0}\,e^{3a(s-u)}\,f(u)\upsilon_\zs{2}(u)\d u$. It is clear that
$$
\vert J_\zs{1,t}(f)\vert \le \check{\varrho}_\zs{2}\Vert f\Vert_\zs{*,n}/3\,.
$$
So,
$$
\left|
a
\int^{t}_\zs{0}
e^{2a(t-s)}\,g(s)
J_\zs{1,s}(f)
\,\d s
\right|\,\le
\check{\varrho}_\zs{2}\Vert f\Vert_\zs{*,t}\,\Vert g\Vert_\zs{*,t}/6\,.
$$
Now we represent  the corresponding integral for $J_\zs{2,t}(f)$  as
$$
\int^{t}_\zs{0}
e^{2a(t-s)}\,g(s)
J_\zs{2,s}(f)
\,\d s
=
\int^{t}_\zs{0}\,e^{3au}\,\wt{J}_\zs{t-u,u}(f,g)\,\d u
$$
and
$$
\wt{J}_\zs{t,u}(f,g)=\int^{t}_\zs{0}\,e^{2a(t-s)}\,g(s+u)\,f(s)\,\upsilon_\zs{2}(s)\d s\,.
$$
Using Lemma \ref{Le.sec:A.7-06-1} and the last inequality in
\eqref{sec:A.8-06-1} we obtain that
$$
\sup_\zs{0\le u\le t}
\left\vert
\wt{J}_\zs{t-u,u}(f)
\right\vert
\le\,
2\varrho^{2}_\zs{*}\,
\Omega_\zs{t}(g,f)
\,,
$$
where
$$
\Omega_\zs{t}(g,f)=\sup_\zs{v\ge 0\,,\,u\ge 0\,,0\le v+u\le t}
\ \
\left\vert
\int^{v}_\zs{0}
g(s+u)f(s)\d s
\right\vert\,.
$$
Therefore,
$$
\left
\vert
a
\int^{t}_\zs{0}
e^{2a(t-s)}\,g(s)
J_\zs{2,s}(f)
\,\d s
\right\vert
\le
\frac{2\varrho^{2}_\zs{*}}{3}\,
\Omega_\zs{t}(g,f)
\le
\frac{2\varrho^{2}_\zs{*}}{3}\,\varpi_\zs{t}(f,g)\,.
$$
Similarly we can get
$$
\left
\vert
a
\int^{t}_\zs{0}
e^{2a(t-s)}\,g(s)
\check{J}_\zs{s}(f)
\,\d s
\right\vert
\le
\frac{4}{3}\,\Omega_\zs{t}(g,\check{\varepsilon}(f))\,.
$$
Note now that for any fixed $v>0$ and $\theta\ge 0t$ with $v+\theta\le t$
we have
$$
\int^{v}_\zs{0}\,g(u+\theta)\,\check{\varepsilon}_\zs{u}(f)\,\d u
=
\frac{a}{2}
\int^{v}_\zs{0}\,e^{as}\,D_\zs{s,\theta}\,\d s
\,,
$$
where
$D_\zs{s,\theta}=\int^{v-s}_\zs{0}\,g(y+\theta+s)\wt{f}(y) \d y$ and
$\wt{f}(y)=f(y)(1+e^{2ay})$.
Integrating par parts yields
\begin{align*}
D_\zs{s,\theta}&=\left(1+e^{2a(t-s)}\right)
\int^{v-s}_\zs{0}g(z+\theta+s)f(z)\d z\\[2mm]
&+2a
\int^{v-s}_\zs{0}\,e^{2au}\,
\int^{u}_\zs{0}\,g(z+\theta+s)\,
f(z)\,
\d z\,
\d u\,.
\end{align*}
This implies, that
$$
\vert D_\zs{s,\theta}(f,g)\vert
\le 3\varpi_\zs{t}(f,g)\,,
$$
i.e. for any $v>0$ and $\theta \ge 0$
\begin{equation}\label{sec:A.9-0}
\left\vert
\int^{v}_\zs{0}\,g(u+\theta)\,\check{\varepsilon}_\zs{u}(f)\,\d u
\right\vert\,
\le
3\varpi_\zs{t}(f,g)/2\,.
\end{equation}
Therefore,
$$
\Omega_\zs{t}(g,\check{\varepsilon}(f))
\le 3\varpi_\zs{t}(f,g)/2\,.
$$
Hence Lemma \ref{Le.sec:A.7-06-2}.
\endproof

\begin{lemma}\label{Le.sec:A.7-06-2-01}
For any mesurable bounded  $[0,+\infty)\to\bbr$ functions $f$ and $g$,
for any $-a_\zs{max}\le a\le 0 $ and for any $t>0$
$$
a^{2}
\left|
\int^{t}_\zs{0}
e^{a(t-s)}\,g(s)
\E\,\wt{I}_\zs{s}(f)\,\wt{I}_\zs{s}(1)
\,\d s
\right|\,\le
4\,
\Vert f
\Vert^{2}_\zs{*,t}
\left(
a_\zs{max}\,\check{\varrho}_\zs{2}
+3\varrho^{2}_\zs{*}
\right)
\,\varpi_\zs{t}(1,g)
\,.
$$
\end{lemma}
\proof
Firstly, note that if $a=0$ then this bound is obvious. Let now $\vert a\vert>0$. Then,
taking into account the representation \eqref{sec:Stc.6+00}
and the bound
$\vert\check{\varepsilon}_\zs{t}(f)\vert\le \Vert f\Vert_\zs{*,t} $
we obtain that
\begin{equation}
\label{sec:Stc.6+01}
\Vert a
A(f)
\Vert_\zs{*,t}
\le \Vert f
\Vert_\zs{*,t}
\left(
a_\zs{max}\,\check{\varrho}_\zs{2}/3
+\varrho^{2}_\zs{*}
\right)\,.
\end{equation}
Thus, from the definition of $\wt{\varkappa}_\zs{u}(f)$ in
\eqref{sec:Stc.7-06-1}
we obtain that
\begin{equation}
\label{sec:Stc.6+02}
\Vert a\wt{\varkappa}(f)
\Vert_\zs{*,t}\,
\le \Vert f
\Vert^{2}_\zs{*,t}
\left(
2a_\zs{max}\,\check{\varrho}_\zs{2}
+6\varrho^{2}_\zs{*}
\right)\,.
\end{equation}

Moreover, note now, that
$$
\int^{t}_\zs{0}
e^{a(t-s)}\,g(s)
\E\,\wt{I}_\zs{s}(f)\,\wt{I}_\zs{s}(1)
\,\d s
= \int^{t}_\zs{0}\,e^{a(t-u)}\,
\wt{\varkappa}_\zs{u}(f)
\,G_\zs{t-u,u}\,\d u\,,
$$
where $G_\zs{T,u}=\int^{T}_\zs{0}\,e^{a z }\,g(z+u)\d z$. The integrating by parts yields
$$
G_\zs{T,u}=\int^{T}_\zs{0}\,g(z+u)\d z+a
\int^{T}_\zs{0}\,e^{a y}
\left(
\int^{y}_\zs{0}\,g(v+u)\d v
\right)
\,\d y
\,.
$$
So, for any $T+u\le t$ we obtain that
$\vert G_\zs{T,u}\vert\,\le 2\varpi_\zs{t}(1,g)$ and, therefore,
$$
\left\vert
\int^{t}_\zs{0}\,e^{a(t-u)}\,
\wt{\varkappa}_\zs{u}(f)
\,G_\zs{t-u,u}\,\d u
\right\vert
\le
\frac{4}{a^{2}}
\Vert f
\Vert^{2}_\zs{*,t}
\left(
a_\zs{max}\,\check{\varrho}_\zs{2}
+3\varrho^{2}_\zs{*}
\right)
\,\varpi_\zs{t}(1,g)\,.
$$
Hence Lemma \ref{Le.sec:A.7-06-2-01}. \endproof

\begin{lemma}\label{Le.sec:A.7-06-3}
For any measurable bounded  $[0,+\infty)\to\bbr$ functions $f$ and $g$,
for any $-\infty< a\le 0 $ and for any $t>0$
\begin{equation}\label{sec:A.9}
\left|
a
\int^{t}_\zs{0}
e^{2a(t-s)}\,
\check{\tau}_\zs{s}(f,g)
\,\d s
\right|\,\le
9\,\varpi^{*}_\zs{t}(f,g)\,.
\end{equation}
\end{lemma}
\proof
Firstly  note, that
using the bound
\eqref{sec:A.9-0}
with $\theta=0$ we obtain
\begin{equation}
\label{sec:A.10}
\left\vert
\int^{t}_\zs{0}g(s)\check{\varepsilon}_\zs{s}(f)\d s
\right\vert
\le 3\varpi_\zs{t}(f,g)/2\,.
\end{equation}
So, $\vert\tau_\zs{t}(f,g)\vert\le 4\varpi_\zs{t}(f,g)$. Moreover, through the bound \eqref{sec:A.10}
and Lemma \ref{Le.sec:A.7-06-1} we obtain that
$$
\left|\,
\int^{t}_\zs{0}
e^{2a(t-s)}\,f(s)
\check{\varepsilon}_\zs{s}(g)
\,\d s
\right|\,
\le\,3\varpi_\zs{t}(f,g)
\,.
$$
Using again  Lemma \ref{Le.sec:A.7-06-1}  and taking into account that $a\tau_\zs{t}(1,1)=(e^{2at}-1)/2$
we estimate
$$
\left|\,a
\int^{t}_\zs{0}
e^{2a(t-s)}\,f(s)
g(s)\tau_\zs{s}(1,1)
\,\d s
\right|\,
\le\,\varpi_\zs{t}(f,g)
\,.
$$
Thus, from taking into account the definition
\eqref{sec:Stc.6-02-0} we obtain the bound \eqref{sec:A.9}.
Hence Lemma \ref{Le.sec:A.7-06-3}.
\endproof

\bigskip

\subsection{Property of Penalty term}

\begin{lemma}\label{Lem.A.1}
For any $n\geq 1$, $\gamma\in\Gamma$ and $0<\varepsilon<1$
\begin{equation}
\label{penalty-00}
P_\zs{n}(\gamma)\leq\frac{\E\,\Er_\zs{n}(\gamma)}{1-\varepsilon}+\frac{\c^{*}_\zs{n}}{n\varepsilon(1-\varepsilon)}
\,.
\end{equation}
\end{lemma}
\proof
By the definition of $\Er_\zs{n}(\gamma)$ one has
\begin{align*}
\Er_\zs{n}(\gamma)&=\sum^{n}_\zs{j=1}\,(\gamma(j)\theta_\zs{j,n}^*-\theta_j)^2
=\sum^{n}_\zs{j=1}\,\left(\gamma(j)(\theta_\zs{j,n}^*-\theta_j)+(\gamma(j)-1)\theta_j\right)^2 \\[2mm]
&
\ge
\sum^{n}_\zs{j=1}\,\gamma(j)^2(\theta_\zs{j,n}^*-\theta_j)^2+
2\sum^{n}_\zs{j=1}\,\gamma(j)(\gamma(j)-1)\theta_j(\theta_\zs{j,n}^*-\theta_j).
\end{align*}
Taking into account the condition $\B_\zs{2})$ and the definition
\eqref{sec:Imp.12} we obtain that the last term in tho sum can be replaced as
$$
\sum^{n}_\zs{j=1}\,\gamma(j)(\gamma(j)-1)\theta_j(\theta_\zs{j,n}^*-\theta_j)
=
\sum^{n}_\zs{j=1}\,\gamma(j)(\gamma(j)-1)\theta_j(\wh{\theta}_\zs{j,n}-\theta_j)
\,,
$$
i.e.
$
\E\,\sum^{n}_\zs{j=1}\,\gamma(j)(\gamma(j)-1)\theta_j(\theta_\zs{j,n}^*-\theta_j)=0$ and, therefore, taking into account the definition
\eqref{sec:Mo.9} we obtain that
\begin{align*}
\E\,\Er_\zs{n}(\gamma)&\geq\sum^{n}_\zs{j=1}\,\gamma(j)^2\E\,(\theta_\zs{j,n}^*-\theta_j)^2=
\sum^{n}_\zs{j=1}\,\gamma(j)^2\E\,\left(\frac{\xi_\zs{j,n}}{\sqrt{n}}-g_\zs{\gamma}(j)\wh{\theta}_\zs{j}\right)^2\\[2mm]
&
\ge
P_\zs{n}(\gamma)-\frac{2}{\sqrt{n}}\E\,\sum^{n}_\zs{j=1}\,\gamma(j)^2g_\zs{\gamma}(j)\wh{\theta}_\zs{j,n}
\xi_\zs{j}
\\[2mm]&
\ge (1-\varepsilon)\,P_{n}(\gamma)
- \frac{1}{\varepsilon}\E\,\sum^{n}_\zs{j=1}\,g_\zs{\gamma}^{2}(j)\wh{\theta}_\zs{j}^{2}
\,.
\end{align*}
The inequality \eqref{sec:Mo.13_Ub++c-n}
 implies the bound \eqref{penalty-00}. Hence  Lemma \ref{Lem.A.1}.
\endproof

\bigskip
\bigskip

\begin{lemma}\label{Le.sec:App.3+1}
For any $m>0$ the term $\R_\zs{0,n}$ introduced in \eqref{sec:Lo.12}
satisfies the following property
\begin{equation}\label{sec:App.5--0}
\lim_\zs{T\to\infty}\,n^{m}\,
\R_\zs{0,n}\,=0\,.
\end{equation}
\end{lemma}
\proof
First, setting $\zeta_\zs{n}=\sum^{d_\zs{n}}_\zs{j=1}\,\kappa^{2}_\zs{j}\,a_\zs{j}$,
we obtain that
$$
\left\{
S_\zs{\kappa}\notin W_\zs{k,\r}
\right\}
=
\left\{
\sum^{d_\zs{n}}_\zs{j=1}\,\kappa^{2}_\zs{j}\sum^{k}_\zs{l=0}\,
\Vert \phi^{(l)}_\zs{j}\Vert^{2}
>\r
\right\}
=
\left\{
\zeta_\zs{n}
>\r
\right\}
\,.
$$
Moreover, note that one can check directly that
$$
\lim_\zs{n\to \infty}\,
\E\,\zeta_\zs{n}=
\lim_\zs{n\to \infty}\,
\frac{1}{v_\zs{n}}
\sum^{d_\zs{n}}_\zs{j=1}\,s^{*}_\zs{j}\,a_\zs{j}=\r_\zs{\ve}=
(1-\ve)\r\,.
$$
So, for sufficiently large $n$ we obtain that
$$
\left\{
S_\zs{\kappa}\notin W_\zs{k,\r}
\right\}
\subset
\left\{
\wt{\zeta}_\zs{n}>
\r_\zs{1}
\right\}
\,,
$$
where $\r_\zs{1}=\r\ve/2$,
$$
\wt{\zeta}_\zs{n}=\zeta_\zs{n}-\E\,\zeta_\zs{n}
=\frac{1}{v_\zs{n}}\,\sum^{d_\zs{n}}_\zs{j=1}\,s^{*}_\zs{j}a_\zs{j}\wt{\eta}_\zs{j}
\quad\mbox{and}\quad
\wt{\eta}_\zs{j}=\eta^{2}_\zs{j}-1\,.
$$
Through the correlation inequality from
\cite{GaPeSPA_2013}
we can get that for any $p\ge 2$ there exists some constant $C_\zs{p}>0$  for which
$$
\E\,\wt{\zeta}^{p}_\zs{n}\le C_\zs{p}
\frac{1}{v^{p}_\zs{n}}\,
\left( \sum^{d}_\zs{j=1}\,
(s^{*}_\zs{j})^{2}a^{2}_\zs{j}
\right)^{p/2}
\le C \,v^{-\frac{p}{4k+2}}_\zs{n}
\,,
$$
i.e.  the expectation
$\E\,\wt{\zeta}^{p}_\zs{n}\to 0$ as $n\to\infty$. Therefore, using the Chebychev inequality
we obtain that for any $m>1$
$$
n^{m}\P(\wt{\zeta}_\zs{n}>\r_\zs{1})\to 0
\quad\mbox{as}\quad
n\to\infty\,.
$$

\noindent
Hence Lemma \ref{Le.sec:App.3+1}. \endproof

\subsection{The van Trees inequality for the Levy processes.}\label{subsec:App.4}

In this section we consider the following continuous time
 parametric regression model
\begin{equation}\label{sec:App.5}
  \d y_t=S(t,\theta)\d t+\d \xi_t\,,
  \quad 0\le t\le n\,,
 \end{equation}
where $\xi_\zs{t}=\varrho_\zs{1}W_\zs{t}+\varrho_\zs{2}z_\zs{t}$
and
$$
S(t,\theta)=
\sum^{d}_\zs{i=1}\,\theta_\zs{i}\,\psi_\zs{i}(t)\,,
$$
with the unknown parameters
$\theta=(\theta_\zs{1},\ldots,\theta_\zs{d})'$. Here we assume that the functions $(\psi)_\zs{1\le j\le d}$ are $1$ periodic and orthogonal functions.

Let us denote by $\nu_\zs{\xi}$ the distribution of the process $(\xi_\zs{t})_\zs{0\le t\le n}$
on the Skorokhod space $\D[0,n]$.  One can check directly that in this space for any parameters $\theta\in\bbr^d$, the distribution
$\P_\zs{\theta}$
of the process \eqref{sec:App.5}
 is absolutely continuous with respect to the
$\nu_\zs{\xi}$
 and the corresponding Radon-Nikodym derivative,
for any function $x=(x_\zs{t})_\zs{0\le t\le T}$ from $\D[0,n]$,
 is defined as
\begin{equation}\label{sec:App.7}
f(x,\theta)=
\frac{\d\P_\zs{\theta}}{\d\nu_\zs{\xi}}(x)=
\exp\left\{\int^{n}_\zs{0}\,\frac{S(t,\theta)}{\varrho_\zs{1}}\,\d x^{c}_\zs{t}
-\,\int^{n}_\zs{0}\,
\frac{S^{2}(t,\theta)}{2\varrho^{2}_\zs{1}}\,
\d t
\right\}
\,,
\end{equation}
where
$$
x^{c}_\zs{t}=
\frac{1}{\varrho_\zs{1}}
\left(
x_\zs{t}
-
\int^{t}_\zs{0}\,\int_\zs{\bbr}\,v\,\left(
\mu_\zs{x}(\d s\,,\d v)
-
\Pi(\d v)\d s
\right)
\right)
$$
and for any measurable set $\Gamma$ in $\bbr$ with $0\notin \Gamma$
$$
\mu_\zs{x}([0,t]\times\Gamma)=\sum_\zs{0\le s\le t}\,
\Chi_\zs{\{\Delta\xi_\zs{s}\in \varrho_\zs{2} \Gamma\}}
\,.
$$

\noindent
Let $\Phi$ be a prior density on $\bbr^d$ having
the following form:
$$
\Phi(\theta)=\Phi(\theta_1,\ldots,\theta_d)=\prod_{j=1}^d\varphi_\zs{j}(\theta_\zs{j})\,,
$$
where $\varphi_\zs{j}$ is some continuously differentiable density in $\bbr$.
Moreover, let $g(\theta)$ be a continuously differentiable $\bbr^d\to \bbr$ function such that,
for each $1\le j\le d$,
\begin{equation}\label{sec:App.8}
\lim_\zs{|\theta_\zs{j}|\to\infty}\,
g(\theta)\,\varphi_\zs{j}(\theta_\zs{j})=0
\quad\mbox{and}\quad
\int_\zs{\bbr^d}\,|g^{\prime}_\zs{j}(\theta)|\,\Phi(\theta)\,\d \theta
<\infty\,,
\end{equation}
where
$$
g^{\prime}_\zs{j}(\theta)=\frac{\partial g(\theta)}{\partial\theta_\zs{j}}\,.
$$
For any $\cB(\cX)\times\cB(\bbr^d)-$
measurable integrable function $H=H(x,\theta)$ we denote
\begin{align*}
\wt{\E}\,H&=\int_{\bbr^d}\,
\int_\zs{\cX}\,H(x,\theta)\,\d \P_\zs{\theta}\,\Phi(\theta) \d \theta\\[2mm]
&=
\int_{\bbr^d}\,\int_\zs{\cX}\,
H(x,\theta)\,f(x,\theta)\,\Phi(\theta)\d \nu_\zs{\xi}(x)\, \d \theta\,,
\end{align*}
where $\cX=\D[0,n]$.

\begin{lemma}\label{Le.sec:App.3}
For any $\cF^y_T$-measurable square integrable function $\wh{g}_\zs{T}$
 and for any $1\le j\le d$, the following inequality holds
$$
\wt{\E}(\wh{g}_\zs{n}-g(\theta))^2\ge
\frac{\Lambda^2_\zs{j}}{n\Vert \psi_\zs{j}\Vert^{2}\varrho^{-2}_\zs{1}+I_\zs{j}}\,,
$$
where
$$
\Lambda_\zs{j}=\int_\zs{\bbr^d}\,g^{\prime}_\zs{j}(\theta)\,\Phi(\theta)\,\d \theta
\quad\mbox{and}\quad
I_\zs{j}=\int_\zs{\bbr}\,\frac{\dot{\varphi}^2_\zs{j}(z)}{\varphi_\zs{j}(z)}\,\d z\,.
$$
\end{lemma}
\noindent {\bf  Proof.}
First of all note that, the density \eqref{sec:App.7}
on the process $\xi$
is bounded
with respect to $\theta_\zs{j}\in\bbr$ and for any $1\le j\le d$
$$
\limsup_\zs{|\theta_\zs{j}|\to\infty}\,f(\xi,\theta)\,=\,0\,.
\quad\quad\mbox{a.s.}
$$
Now, we set
$$
\wt{\Phi}_\zs{j}=\wt{\Phi}_\zs{j}(x,\theta)=
\frac{\partial\,(f(x,\theta)\Phi(\theta))/\partial\theta_\zs{j}}{f(x,\theta)\Phi(\theta)}
 \,.
$$
Taking into account the condition \eqref{sec:App.8} and
integrating by parts yield
\begin{align*}
\wt{\E}&\left((\wh{g}_\zs{n}-g(\theta))\wt{\Phi}_\zs{j}\right)
=\int_{\cX\times\bbr^d}\,(\wh{g}_\zs{n}(x)-g(\theta))\frac{\partial}{\partial\theta_\zs{j}}
\left(f(x,\theta)\Phi(\theta)\right)\d \theta\,\nu_\zs{\xi}(\d x)\\[2mm]
&=\int_{\cX\times\bbr^{d-1}}\left(\int_{\bbr}\,
g^{\prime}_\zs{j}(\theta)\,
f(x,\theta)\Phi(\theta)\d \theta_\zs{j}\right)\left(\prod_{i\neq j}\d \theta_i\right)\,\nu_\zs{\xi}(\d x)
=\Lambda_\zs{j}\,.
\end{align*}
Now by the Bouniakovskii-Cauchy-Schwarz inequality
we obtain the following lower bound for the quadratic risk
$$
\wt{\E}(\wh{g}_\zs{T}-g(\theta))^2\ge
\frac{\Lambda^2_\zs{j}}{\wt{\E}\Psi_\zs{j}^2}\,.
$$
To study the denominator in the left hand of this
inequality note that in view of the representation
 \eqref{sec:App.7}
$$
\frac{1}{f(y,\theta)}
\frac{\partial\,f(y,\theta)}{\partial\theta_\zs{j}}
=\frac{1}{\varrho_\zs{1}}\,
\int^{n}_\zs{0}\,\psi_\zs{j}(t)\,\d w_\zs{t}\,.
$$
Therefore, for each $\theta\in\bbr^d$,
$$
\E_\zs{\theta}\,
\frac{1}{f(y,\theta)}
\frac{\partial\,f(y,\theta)}{\partial\theta_\zs{j}}
\,
=0
$$
and
$$
\E_\zs{\theta}\,
\left(
\frac{1}{f(y,\theta)}
\frac{\partial\,f(y,\theta)}{\partial\theta_\zs{j}}
\right)^2
=\,
\frac{1}{\varrho^{2}_\zs{1}}
\int^{n}_\zs{0}\,\psi^2_\zs{j}(t)\d t
=
\frac{n}{\varrho^{2}_\zs{1}}
\Vert\psi\Vert^{2}
\,.
$$
Taking into account that
$$
\wt{\Phi}_\zs{j}=
\frac{1}{f(x,\theta)}
\frac{\partial\,f(x,\theta)}{\partial\theta_\zs{j}}
+
\frac{1}{\Phi(\theta)}
\frac{\partial\,\Phi(\theta))}{\partial\theta_\zs{j}}
 \,,
$$
we get
$$
\wt{\E}\Psi_\zs{j}^2=
\frac{n}{\varrho^{2}_\zs{1}}\,\Vert\psi\Vert^{2}
+\,I_\zs{j}\,.
$$
Hence
Lemma~\ref{Le.sec:App.3}.
\endproof

\end{document}